\documentclass[preprint,5p]{elsarticle}

\pdfoutput=0

\usepackage{amsmath,amsfonts,amssymb}
\usepackage{times}
\usepackage[english]{babel}
\usepackage{bm}
\usepackage{algorithm}
\usepackage{algorithmic}
\usepackage{url}
\usepackage{multirow}

\def\XS{\xspace}

\RequirePackage{ifthen}
\RequirePackage{calc}
\RequirePackage{xspace}

\def\<{\og}
\def\>{\fg}

\newcommand{\taille}[1][\scad]{%
\ifthenelse{#1 = -5}{}{}%
\ifthenelse{#1 = -4}{\tiny}{}%
\ifthenelse{#1 = -3}{\scriptsize}{}%
\ifthenelse{#1 = -2}{\footnotesize}{}%
\ifthenelse{#1 = -1}{\small}{}%
\ifthenelse{#1 = 0}{\normalsize}{}%
\ifthenelse{#1 = 1}{\large}{}%
\ifthenelse{#1 = 2}{\Large}{}%
\ifthenelse{#1 = 3}{\LARGE}{}%
\ifthenelse{#1 = 4}{\huge}{}%
\ifthenelse{#1 = 5}{\Huge}{}}

\newboolean{@serif}
\setboolean{@serif}{true}
\def\scad{-5} 

\newcounter{taille}
\newcommand{\sca}[2][\scad]{\setcounter{taille}{#1}%
\ifthenelse{\boolean{@serif}}
{{\taille[\thetaille]\textsc{#2}}}
{\setcounter{taille}{\value{taille}-1}{\uppercase{\taille[\thetaille]#2}}}}

\def\rem#1{}                    

\def\reff#1{\XS(\ref{#1})}

\def\BibTeX{{\rm B\kern-.05em{\sc i\kern-.025em b}\kern-.08em
T\kern-.1667em\lower.7ex\hbox{E}\kern-.125emX}}

\newcommand{\beaunom}[4][0]{{\taille[#1]#2 \uppercase{#3}}%
{\sca[#1]{#4}}\XS}

\newcommand{\cnrs}[1][\scad]{\sca[#1]{cnrs}\XS}

\def\UPS{Universit\'e de Paris--Sud\XS}
\newcommand{\ups}[1][\scad]{\sca[#1]{ups}\XS}

\newcommand{\upvi}[1][\scad]{\sca[#1]{upvi}\XS}

\newcommand{\drm}[1][\scad]{\sca[#1]{drm}\XS}

\newcommand{\dsv}[1][\scad]{\sca[#1]{dsv}\XS}

\newcommand{\shfj}[1][\scad]{\sca[#1]{shfj}\XS}

\newcommand{\lena}[1][\scad]{\sca[#1]{lena}\XS}
\newcommand{\sxqt}[1][\scad]{\sca[#1]{{\normalsize 640}}\XS}

\newcommand{\ESE}[1][\scad]{\sca[#1]{sup\'elec}\XS}

\newcommand{\adresscea}[1][\scad]{4 Place du G\'en\'eral Leclerc, 91406 Orsay Cedex, France\XS}
\def\helvnormal#1{{\fontfamily{phv}\fontsize{11}{14pt}\selectfont #1}}
\def\helvfoot  #1{{\fontfamily{phv}\fontsize{9}{14pt} \selectfont #1}}
\def\helvscript#1{{\fontfamily{phv}\fontsize{8}{14pt} \selectfont #1}}
\def\pth#1{\left(#1\right)}                
\def\acc#1{\left\{#1\right\}}              
\def\cro#1{\left[#1\right]}                
\def\bars#1{\left|#1\right|}

\def\bigpth#1{\bigl(#1\bigr)}              
\def\bigacc#1{\bigl\{#1\bigr\}}            
\def\bigcro#1{\bigl[#1\bigr]}

\def\cond{\textrm{Cond}\,}

\def\rang{{\mathrm{rang}}\,}                       
\def\ker{\textrm{Ker}\,}                         
\def\img{\textrm{Im}\,}                          
\def\vect{{\mathrm{Vect}}\,}                       
\def\sgn{{\mathrm{sgn}}}                         

\def\card{{\mathrm{Card}}}      

\def\WHERE{\text{where\:}}

\newsavebox{\fminibox}
\newlength{\fminilength}
\newenvironment{fminipage}[1][\linewidth]
{\setlength{\fminilength}{#1}
\begin{lrbox}{\fminibox}\begin{minipage}{\fminilength}}
{\end{minipage}\end{lrbox}\noindent\fbox{\usebox{\fminibox}}}

\def\M{^{-1}} \def\T{^\tD} \def\+{^\dagger}
\def\I{\,|\,}           

\def\ldotsv{,\,\ldots,\,}

\def\nequiv{\not\kern-.05em\equiv}
\def\egal{\kern-.5em=\kern-.5em}        
\def\propt{\kern-.2em\propto\kern-.2em} 
\def\wh#1{\widehat{#1}}                 
\def\wt#1{\widetilde{#1}} 

\def\froc#1#2{{#1/#2}}                  

\def\intdouble{\int\kern-0.3em\int}
\def\inttriple{\int\kern-0.3em\int\kern-0.3em\int}

\def\rond#1{\overset{\kern-0.33em~_\circ}{#1}}
\def\rondit[#1]#2{\overset{\kern#1~_\circ}{#2}}

\def\incirc#1{\pscirclebox[framesep=1pt]{\scriptsize#1}}

\def\sss{\scriptscriptstyle}

\RequirePackage{amsmath}
\RequirePackage{xspace}
\RequirePackage{bbm}

\def\XS{\xspace}
\DeclareMathAlphabet{\mathb}{OML}{cmm}{b}{it}
\def\sbm#1{\ensuremath{\mathb{#1}}}                
\def\sbmm#1{\ensuremath{\boldsymbol{#1}}}          
\def\sdm#1{\ensuremath{\mathrm{#1}}}               
\def\sbv#1{\ensuremath{\mathbf{#1}}}               
\def\scu#1{\ensuremath{\mathcal{#1\XS}}}           
\def\scb#1{\ensuremath{\boldsymbol{\mathcal{#1}}}} 
\def\sbl#1{\ensuremath{\mathbbm{#1}}}              

  \def\ab{{\sbm{a}}\XS}

\def\Eb{{\sbm{E}}\XS}  
\def\Fb{{\sbm{F}}\XS}  \def\fb{{\sbm{f}}\XS}

\def\Lb{{\sbm{L}}\XS}  
  
  \def\nb{{\sbm{n}}\XS}

  \def\tb{{\sbm{t}}\XS}
\def\Ub{{\sbm{U}}\XS}  \def\ub{{\sbm{u}}\XS}
\def\Vb{{\sbm{V}}\XS}  \def\vb{{\sbm{v}}\XS}
\def\Wb{{\sbm{W}}\XS}  
  
\def\Yb{{\sbm{Y}}\XS}

\def\Ac{{\scu{A}}\XS}

\def\Dc{{\scu{D}}\XS}   
   
\def\Fc{{\scu{F}}\XS}

\def\Jc{{\scu{J}}\XS}

\def\Mc{{\scu{M}}\XS}   
\def\Nc{{\scu{N}}\XS}

\def\Tc{{\scu{T}}\XS}   
\def\Uc{{\scu{U}}\XS}   
   
\def\Wc{{\scu{W}}\XS}

  \def\tD{{\sdm{t}}\XS}

\def\Cbb{{\sbl{C}}\XS}  
  
\def\Ebb{{\sbl{E}}\XS}

\def\Nbb{{\sbl{N}}\XS}

\def\Qbb{{\sbl{Q}}\XS}  
\def\Rbb{{\sbl{R}}\XS}

\def\Zbb{{\sbl{Z}}\XS}

\def\babs{\begin{abstract}}             \def\eabs{\end{abstract}}
\def\barr{\begin{array}}                \def\earr{\end{array}}
\def\bcc{\begin{center}}                \def\ecc{\end{center}}
\def\cl#1{\centerline{#1}}
\def\bdes{\begin{description}}          \def\edes{\end{description}}
\def\bdoc{\begin{document}}             \def\edoc{\end{document}}
\def\ben{\begin{enumerate}}             \def\een{\end{enumerate}}
\def\beqn{\begin{eqnarray}}             \def\eeqn{\end{eqnarray}}
\def\beqnl#1{\beqn\label{#1}}           \def\eeqnl#1{\label{#1}\eeqn}
\def\beqnx{\begin{eqnarray*}}           \def\eeqnx{\end{eqnarray*}}
\def\bseqn{\begin{subeqnarray}}         \def\eseqn{\end{subeqnarray}}
\def\beq#1\eeq{\begin{equation}#1\end{equation}}
\def\bal#1\eal{\begin{align}#1\end{align}}
\def\balx#1\ealx{\begin{align*}#1\end{align*}}
\def\beqx{$$}                           \def\eeqx{$$}
\def\bfig{\protect\begin{figure}}       \def\efig{\protect\end{figure}}
\def\bfigx{\protect\begin{figure*}}     \def\efigx{\protect\end{figure*}}
\def\bfigt{\protect\begin{figurette}}   \def\efigt{\protect\end{figurette}}
\def\bfl{\begin{flushleft}}             \def\efl{\end{flushleft}}
\def\bfr{\begin{flushright}}            \def\efr{\end{flushright}}
\def\bit{\begin{itemize}}               \def\eit{\end{itemize}}
\def\bmi{\begin{minipage}}              \def\emi{\end{minipage}}
\def\bfmi{\begin{fminipage}}            \def\efmi{\end{fminipage}}
\def\bpic{\begin{picture}}              \def\epic{\end{picture}}
\def\bqu{\begin{quote}}                 \def\equ{\end{quote}}
\def\bqun{\begin{quotation}}            \def\equn{\end{quotation}}
\def\bsl{\begin{slide}}                 \def\esl{\end{slide}}
\def\btabb{\begin{tabbing}}             \def\etabb{\end{tabbing}}
\def\btabl{\begin{table}}               \def\etabl{\end{table}}
\def\btablx{\begin{table*}}             \def\etablx{\end{table*}}
\def\btab{\begin{tabular}} 
\def\btabu{\begin{tabular}}             \def\etabu{\end{tabular}}
\def\btabx{\begin{tabular*}}            \def\etabx{\end{tabular*}}
\def\bbib{\begin{thebibliography}{999}} \def\ebib{\end{thebibliography}}
\def\bver{\begin{verbatim}}             \def\ever{\end{verbatim}}
\def\bca{\begin{cases}}                          \def\eca{\end{cases}}

\def \E {{\mathbb{E}}}
\def \R {\mathbb{R}}

\def\cc#1{\setlength{\tabcolsep}{0pt}\btabu{c}#1\etabu}
\newcommand{\figc}[2][]
{\setlength{\tabcolsep}{0pt}\btabu{c}\includegraphics[#1]{#2}\etabu}
\def\yaxis#1{\cc{\rotatebox{90}{{\small #1}}}}

\def\cl#1{\centerline{#1}}
\def\eqcl#1{\cl{$\displaystyle#1$}}

\usepackage{xcolor}
\usepackage[normalem]{ulem}
\newcommand{\ADDED}[1]{{#1}}
\newcommand{\DELETED}[1]{\textcolor{red}{\sout{#1}}}
\newcommand{\REPLACE}[2]{\textcolor{red}{\sout{}}\textcolor{blue}{{#2}}}

\graphicspath{{./figures_final/}}

\newenvironment{statement}[1]{\begin{description}
\item[\bfseries #1:]}{\end{description}}
\journal{Fontiers in Fractal Physiology}

\newlength{\largEncadre}
\newcommand{\encFunc}[1]{ 
\settowidth{\largEncadre}{#1}
\fcolorbox{red}{red!20}{\parbox{\largEncadre}{#1}}}
\newcommand{\encArt}[1]{ 
\settowidth{\largEncadre}{#1}
\fcolorbox{blue}{blue!20}{\parbox{\largEncadre}{#1}}}
\newcommand{\encFuz}[1]{ 
\settowidth{\largEncadre}{#1}
\fcolorbox{green}{green!20}{\parbox{\largEncadre}{#1}}}

\begin{document}

\title{Scale-free and multifractal time dynamics of fMRI signals during rest and task}

\author[cea]{P.~Ciuciu\corref{corresponding}}
\author[cea,parietal,unicog]{G.~Varoquaux}
\author[enslyon]{P.~Abry}
\author[desposito]{S.~Sadaghiani}
\author[unicog]{A.~Kleinschmidt}

\cortext[corresponding]{Corresponding author (philippe.ciuciu@cea.fr).}
\address[cea]{CEA/Neurospin, B\^at 145, F-91191 Gif-Sur-Yvette.}
\address[parietal]{Parietal project-team, INRIA Saclay-\^ile de France, F-91191 Gif-Sur-Yvette.}
\address[unicog]{INSERM U992, NeuroSpin center, F-91191 Gif-sur-Yvette.}
\address[enslyon]{Physics Dept. (CNRS, UMR 5672) ENS Lyon, 46 all\'ee d'Italie, F-69007 Lyon.}
\address[desposito]{D'esposito Lab, University of California, 132 Barker Hall Berkeley CA-94720.}

\begin{abstract}

Scaling temporal dynamics in functional MRI~(fMRI) signals have been evidenced for a decade as intrinsic
characteristics of ongoing brain activity~\cite{Zarahn97}. 
Recently, scaling properties were shown to fluctuate across brain networks and to be modulated between rest
and task~\cite{He11}: 
Notably, Hurst exponent, quantifying long memory, decreases under task in activating and deactivating brain
regions. In most cases, such results were obtained: First, from univariate~(voxelwise or regionwise) analysis,
hence focusing on specific cognitive systems such as Resting-State
Networks~(RSNs) and raising the issue of the specificity of this scale-free dynamics modulation in RSNs.
Second, using analysis tools designed to measure a single scaling exponent related to the second order statistics of
the data, thus relying on models that either implicitly or explicitly assume Gaussianity and~(asymptotic)
self-similarity, while fMRI signals may significantly depart from those either of those two
assumptions~\cite{Ciuciu08,Wink08}.

To address these issues, the present contribution elaborates on the analysis of the scaling properties of fMRI
temporal dynamics by proposing two significant variations.
First, scaling properties are technically investigated using the recently introduced Wavelet
Leader-based Multifractal formalism (WLMF)~\cite{WENDT:2007:E,ABRY:2012:A}. 
This measures a collection of scaling exponents, thus enables a richer and more versatile description of scale
invariance (beyond correlation and Gaussianity), referred to as multifractality.
Also, it benefits from improved estimation performance compared to tools previously used in the literature.
Second, scaling properties are investigated in both RSN and non-RSN structures~(e.g.,  artifacts), at a broader spatial
scale than the voxel one, using a multivariate approach, namely the Multi-Subject Dictionary Learning~(MSDL)
algorithm~\cite{Varoquaux11} that produces a set of spatial components that appear
more sparse than their Independent Component Analysis~(ICA) counterpart.

These tools  are combined and applied to a fMRI dataset comprising 12 subjects with resting-state and activation
runs~\cite{Sadaghiani09}. 
Results stemming from those analysis confirm the already reported task-related decrease of long memory in functional
networks,  but also show that it occurs in artifacts, thus making this feature not specific to functional networks.
Further, results indicate that most fMRI signals appear multifractal at rest except in non-cortical regions.
Task-related modulation of multifractality appears only significant
in functional networks and thus can be considered as the key property disentangling functional
networks from artifacts. 
These finding are discussed in the light of the recent literature reporting scaling dynamics of EEG microstate
sequences at rest and addressing non-stationarity issues in temporally independent fMRI modes.

\end{abstract}


\maketitle

\sloppy 
\section{Introduction}\label{sec:intro}

Much of what is known about brain function stems from studies in which a task or a stimulus is
administred and the resulting changes in neuronal activity and behaviour are measured.
From the advent of human electroencephalography~(EEG) to cognitive activation paradigms in
functional Magnetic Resonance Imaging~(fMRI), this approach proved very successful to study brain function, and more
precisely functional specialization in human brain.
It has relied, on one hand, on contrasting signal magnitude between different experimental conditions~\cite{Rosen98} or
task-specific hemodynamic response~(HRF) shape~\cite{Dale99} and, on other-hand, on statistical me\-thods often
framed within linear or bilinear modelling strategies~\cite{Friston95e,Makni05,Makni08}.

Spontaneous modulations of neural activity in Blood Oxygenation Level Dependent~(BOLD) fMRI signals however arise
without external input or stimulus and thus depict intrinsic brain activity~\cite{Damoiseaux06}.
This ongoing activity constitutes a major part of fMRI recordings and is responsible for most of
brain energy consumption. 
It has hence been intensively studied over the last decade using various methods ranging from
\emph{univariate}, i.e., Seed-based linear Correlation Analysis~(SCA)~\cite{Biswal95,Greicius03}, 
to \emph{multivariate} methods such as Independent Component Analysis~(ICA)~\cite{Calhoun01,Beckmann04},
group-level ICA~\cite{Cole10,Varoquaux10} or more recent dictionary learning techniques~\cite{Varoquaux11}.
All these methods have revealed that interactions between brain regions, also referred to as
\emph{functional connectivity},
occur through these spontaneous modulations and consistantly vary between rest and task~\cite{Damoiseaux06,Fox07}.
Resting-State Network~(RSN) extraction from resting-state fMRI time series is thus achieved either by thresholding
the correlation matrix computed between voxels or regions~(seed-based or univariate approach) or by identifying
spatial maps in ICA-based algorithms that closely match RSNs such as somato-sensory systems~(visual, motor, auditory),
the default mode and attentional networks~(ventral and dorsal)~\cite{Fox07,Smith09}.
For a recent review about the pros and cons of the
SCA and ICA approaches to RSN extraction, the reader can refer to~\cite{Cole10}. Once RSNs are extracted,
their topological properties can be analyzed with respect to \emph{small-world} or
\emph{scale-free} models~\cite{Chialvo04,Eguiluz05,Zemanova06,Bullmore09}.

In parallel and alternatively to brain topology, the temporal dynamics of brain activity have also been extensively studied.
It is now well accepted that brain activity, irrespective of the imaging technique involved in observation, is always
arrhythmic and shows a scaling, or scale invariant or scale-free, time dynamics, which implies that no time scale plays
a predominant or specific role. Often, scale invariance or \emph{scale-free dynamics} is associated with
\emph{long range correlation} in \emph{time}~\cite{Linkenkaer-Hansen01,Stam04,VandeVille10}, 
and accordingly, in the frequency domain, related to a \emph{power-law} decrease of the power
spectrum~($\Gamma(f)\!\!\propto\!\! 1/f^\beta$ with $\beta\!>\!0$) in the limit of small frequencies~($f\!\rightarrow\! 0$). 
Interestingly, it is generally admitted that only low frequencies~($<0.1$Hz) convey information related to neural
connectivity in fMRI signals~\cite{Cordes01,Leopold03,Achard06}.
Evidence of fractal or \emph{scale-free} behavior in fMRI signals has been
demonstrated for a long while~\cite{Zarahn97,Bullmore01,Bullock03} though it was initially regarded as noise. 
Deeper investigations of the temporal \emph{scale-free} property in fMRI have demonstrated that this
constitutes an intrinsic feature of ongoing brain activity~(cf. e.g.~\cite{Thurner03,Shimizu04,Maxim05,Ciuciu08,Wink08,He10,He11}).
First attempts to identify stimulus-induced signal changes from scaling parameters were proposed
in~\cite{Thurner03,Shimizu04},
where a voxel-based fluctuation analysis was applied to high temporal resolution fMRI data. 
Interestingly, fractal features of voxel time series have enabled to discriminate white matter, 
cerebrospinal fluid and active from inactive brain regions during a block paradigm~\cite{Shimizu04}.
Further, it was shown that scaling properties can be modulated in neurological
disorder~\cite{Maxim05} or between rest and task~\cite{Thurner03,Shimizu04,Ciuciu08,Wink08,He11}: It was shown that
long memory, as quantified by the Hurst exponent, decreases during task in activating and deactivating
brain regions. 
Analyzing scale invariance in temporal dynamics may thus provide new insights into how the brain works by mapping
quantitative estimations of parameters with good specificities to cognitive states, task
performance~\cite{Shimizu04,Wink08,He10,He11}.

Small world and scale-free topology led to model brain as a complex critical system, that is as a large conglomerate of
interacting components, with possibly nonlinear interactions~\cite{Bak95,Chialvo10}.
Further, these complex systems were then regarded as potential origins for long-range correlation spatio-temporal patterns,
as critical systems, i.e., complex systems driven close to their phase transitions, constitute known mechanism yielding
scaling time dynamics and generic $1/f$ power spectral densities~(see e.g.~\cite{Chialvo10}). 
They however so far failed to account for the existence of possibly richer scaling properties (such as e.g.,
multifractality). At a general level, scale invariance in time dynamics and scale-free property of brain topology are,
in essence, totally independent properties that must not be confused one with the other.
Whether or not and how these two scale-free instances are related one to the other in the fMRI context remains a difficult
and largely unsolved issue, far beyond the scope of the present contribution, that concentrates instead on performing
a thorough analysis of scale invariance temporal dynamics in fMRI signals. 

In the existing literature, the analysis of scale invariance in fMRI signals suffers from two limitations: 
First, it has often been performed at the voxel or region level,
thus consisting of a collection of univariate analyses, suffering from the classical bias of voxel selection or region
definition. Moreover, although the fluctuation of scale-free dynamics with tissue type has been studied
in~\cite{Shimizu04,Wink08} to derive that stronger persistency occurs in grey matter and that this background activity
might represent neuronal dynamics, no systematic analysis has been undertaken to disentangle the scale-free properties
of RSN and non-RSN components, such as artifacts. This investigation can be better handled using multivariate or ICA-like
approaches.
Second, scale invariance in fMRI signals has mostly been based on spectral analysis and/or
\emph{Detrended Fluctuation Analysis}~(cf. e.g.~\cite{Thurner03,Stam04,He11}).
This amounts to considering that scaling is associated only with the correlation or the spectrum (hence with the second order
statistics) of the data and thus, implicitly and sometimes even explicitly, to assuming Gaussianity and (asymptotic)
self-similarity for the data (cf. e.g.~\cite{eke2002} for a survey in the fMRI context). 
Also, it is now well-known that such technics lack robustness to disentangle stationarity/non-stationarity versus
true scaling property issues and do not allow simple extension to account for richer scaling properties such as those
observed in multifractal models. 
It is well-accepted that wavelet analysis based analysis of scaling (cf.
e.g.~\cite{agf95,avf98,va01,Bullmore01,Fadili02}) yield not only better estimation performance, but also show
significant practical robustness, notably to non-stationarity, while paving the way toward the analysis of scaling
properties beyond the strict second-order (hence beyond Gaussianity and asymptotic self-similarity). 

In this context, the present contribution elaborates on earlier works dedicated to the analysis of scale invariance
in fMRI temporal dynamics by proposing two significant variations. 

First, scale invariance dynamics is not investigated at the voxel or region spatial scale level independently. 
Instead, group-level resting-state networks are segmented by an exploratory multivariate decomposition approach,
namely the MSDL algorithm~\cite{Varoquaux11}, detailed in Section~\ref{sec:msdl}: It produces both
a set of spatial components and a set of times series, for each component and each subject,
that conveys ongoing dynamics in functional networks but also in artifacts. As shown in~\cite{Varoquaux11},
the sparsisty promoting regularization involved in the MSDL algorithm enables to recover less noisy spatial maps
than group-level or canonical ICA~\cite{Varoquaux10}. This makes their interpretation easier in the context of small
group of individuals. This technique is detailed in Section~\ref{sec:msdl}.

Second, to enable an in-depth analysis of the scaling properties of the temporal dynamics in fMRI signals, 
we resort to multifractal analysis, that measures not a single but a collection of scaling exponents, thus enabling a richer and more versatile description of scale invariance (beyond correlation and Gaussianity), referred to as multifractality.
It is thus likely to better account for the variety and complexity of potential scaling dynamics, as already suggested in
the context of fMRI in e.g.~\cite{Ciuciu08,Wink08}.
However, in contrast to~\cite{Wink08}, and following the track opened in~\cite{Ciuciu08}, we use a recent statistical
analysis tool, the Wavelet Leader-based Multifractal formalism (WLMF)~\cite{WENDT:2007:E,ABRY:2012:A}.
This formalism benefits from better mathematical grounding and shows improved estimation performance compared to tools
previously used in the literature. 
This framework is introduced in Section~\ref{sec:scaling}, after a review of the intuition, models and methodologies
underlying the definition and analysis of scaling temporal dynamics, thus, to some extend, continuing and renewing
the surveys provided in \cite{eke2002,Ciuciu08}.

These tools are combined together and applied to two datasets, corresponding to resting-state and activation runs.
They are described in Section~\ref{sec:data} (see also~\cite{Sadaghiani09}).
Modulations of scale-free and multifractal properties in space, i.e., between functional and artifactual components
but also between rest and task, are statistically assessed at the group level in Section~\ref{sec:results}.

In agreement with findings in~\cite{He11}, the results reported here confirm that fMRI
signals can be modeled as stationary processes, as well as the decrease of the estimated long memory parameter under task.
However, this is found to occur everywhere in the brain and not specifically in functional networks. 
Moreover, evidence for multifractality in resting-state fMRI signals is demonstrated except for non-cortical regions. 
Task-related modulations of multifractality appear only significant
in functional networks and thus become the key property to disentangle functional
networks from artefacts. 
However, in contrast to what happens for the long memory parameter, this modulation is not monotonous across the brain
and varies between cortical and non-cortical regions.
These results are further discussed in Section~\ref{sec:discussion} in the light of recent findings
related to scale-free dynamics of EEG microstate sequences and non-stationarity of functional modes. 
Conclusions are drawn in Section~\ref{sec:conclusion}.

\section{Data acquisition and analysis}\label{sec:data}

\subsection{Data acquisition}

Twelve right-handed normal-hearing subjects (two female; ages, 19--30) gave written informed consent
before participation in an imaging study on a 3T MRI whole-body scanner (Tim-Trio; Siemens).
The study received ethics committee approval by the authorities responsible for our institution.
Anatomical imaging used a T1-weighted magnetization-prepared rapid acquisition gradient echo sequence
[176 slices, repetition time~(TR) 2300 ms, echo time~(TE) 4.18 ms, field of view~(FOV) 256,
voxel size $1\!\!\times\!\! 1\!\!\times\!\! 1$mm$^3$).
Functional imaging used a $\text{T2}^*$-weighted gradient-echo, echo-planar-imaging sequence~(25 slices,
$\text{TR} \!\!=\!\! 1500$~ms, $\text{TE}\!\!=\!\!30$~ms, FOV 192, voxel size $3\!\!\times\!\! 3\!\!\times\!\! 3$mm$^3$).
Stimulus presentation and response recording used the Cogent Toolbox~(John Romaya, Vision Lab,
UCL\footnote{\url{www.vislab.ucl.ac.uk}}) for Matlab and sound delivery a commercially available MR-compatible
system~(MR Confon).

The rs-fMRI dataset we consider in this study has already been published in~\cite{Sadaghiani09}.
820 volumes of task-free ``resting state'' data~(with closed, blind-folded eyes) were acquired before
getting experimental runs of 820 volumes each.
These experimental runs, which have not been analyzed in~\cite{Sadaghiani09},
involve an auditory detection task~(run 2, motor response), and make
use of a sparse \emph{supra}-threshold auditory stimulus detection.

The auditory stimulus was a 500 ms noise burst with its frequency band modulated at 2 Hz~(from white noise to
a narrower band of 0--5~kHz and back to white noise).
Inter-stimulus intervals ranged unpredictably from 20 to 40 s, with each specific interval used only once.
Subjects were instructed to report as quickly and accurately as possible by a right-hand key
press whenever they heard the target sound despite scanner's background noise.
Details about the definition of each subject's auditory threshold are available in~\cite{Sadaghiani09}.

\subsection{Data analysis}

We used here statistical parametric mapping (SPM5, Wellcome Department of Imaging
Neuroscience, UK\footnote{\url{ww.fil.ion.ucl.ac.uk})}.
for image preprocessing~(realignment, coregistration, normalization to MNI stereotactic space, spatial smoothing
with a 5~mm full-width at half-maximum isotropic Gaussian kernel for single-subject and group analyses)
and our own software developments for subsequent analyses. More precisely, the MSDL algorithm relies on
the \texttt{scikit-learn} Python
toolbox\footnote{(\url{http://scikit-learn.org/stable/}).} and the multifractal analysis
on the \texttt{WLBMF} Matlab toolbox\footnote{(\url{http://perso.ens-lyon.fr/herwig.wendt/})}.

\section{Multivariate decomposition of resting state networks}\label{sec:msdl}

\subsection{Multisubject spatial decomposition techniques}

The fMRI signal observed in a voxel reflects many different processes, such as
cardiac or respiratory noise, movement effects, scanner artifacts, or the BOLD
effect that reveals the underlying neural activity of interest. We separate
these different contributions making use of a recently introduced multivariate analysis technique that
estimates jointly spatial maps and time series characteristic of these different
processes~\cite{Varoquaux11}. Formally, this estimation procedure amounts
to finding $K$ spatial maps $\Vb_s \in \mathbb{R}^{p \times K}$ and
the corresponding time series $\Ub_s \in \mathbb{R}^{n \times K}$,  whose linear
combination fits well the observed brain signals, $\Yb_s \in \mathbb{R}^{n\times p}$, 
of length $n$, measured over $p$ voxels, for subject $s$: 
\begin{align}
    \qquad \Yb_s &= \Ub_s \Vb_s\T + \Eb_s, 
    \label{eq:subject_level}
\end{align}
with $\Eb_s \in \mathbb{R}^{n\times p}$ the subject-level noise, or
residuals not explained by the model. Finding $\Vb_s\T$ enables the
separation of the contributions of the different process that are
mixed at the voxel level, but implies to work on spatial maps
rather than on specific voxels. 
The number of spatial maps, $K$, is not chosen a priori, but selected by the procedure.

This problem can be seen as a blind source separation task in the
presence of noise, and has often been tackled in fMRI using ICA,
combined with principal component analysis (PCA) to reject
noise~\cite{mckeown1998,kiviniemi2003,Beckmann04}. In
the multi-subject configuration, estimating the spatial maps on all
subjects simultaneously makes it easy to relate the factors estimated
across the different subjects. This can be done by concatenating the data
across subject, modeling a common distribution~\cite{Calhoun01}, or by
extending the data-reduction step performed in the PCA by a second level
capturing inter-subject variability~\cite{Varoquaux10}. More recently,
it was proposed that the key to the success of ICA on fMRI data, is to
recover sparse spatial maps~\cite{Daubechies09,Varoquaux10b}. This
hypothesis can be formulated as a sparse prior in model~\eqref{eq:subject_level},
which can then be estimated using sparse PCA or
sparse dictionary learning procedures. With regards to our goal in this
study, extracting time-series specific to the various processes observed,
a strong benefit of such procedures is that they can perform data
reduction, \emph{i.e.}, estimation of the residuals not explained by the
model, and extraction of the relevant signals in a single step informed
by our prior. On the opposite, with ICA-based procedures, the residuals
are selected by the PCA step, and not the ICA step.

\subsection{Multi-subject Dictionary Learning algorithm}

In addition,~\citet{Varoquaux11} have adapted the dictionary learning
procedures to a multi-subject setting, in a so-called \emph{multi-subject
dictionary learning}~(MSDL) framework. On fMRI datasets, the procedure
extracts a group-level atlas of spatial signatures of the processes
observed, as well as corresponding subject-level maps, accounting for the
individual specificities. They show that, with a small spatial smoothness
prior added to the sparsity prior on the maps, the extracted patterns
correspond to the segmentation of various structures in the signal:
functional regions, blood vessels, interstitial spaces, sub-cortical
structures... In these settings, the subject-level maps ${\Vb_s}$ are
modeled as generated by group-level maps $\Vb\in
\mathbb{R}^{p\times K}$ with additional inter-subject variability
that appears as residual terms, $\Fb_s \in \mathbb{R}^{p\times K}$,
at the group level:\\
\eqcl{\forall s \in \acc{1\ldotsv S}, \:\Vb_s = \Vb + \Fb_s.}

The model is estimated by finding the group-level and subject-level maps
that maximize the probability of observing the data at hand with the
given prior. This procedure is known as a Maximum A Posteriori~(MAP) 
estimate, and boils down to minimizing the negated log-likelihood of the model 
with an additional penalizing term. If the two 
sources of unexplained signal, \emph{i.e.} subject-level residuals 
$\Eb_s$ and inter-subject variability $\Fb_s$ are modeled as Gaussian 
random variates, the log-likelihood term is the sum of squares of these
errors. The prior term appears as the sum of the sparsity-inducing 
$\ell_1$ norm of $\Vb$, and the $\ell_2$ norm of the gradient of the
map, enforcing the smoothness. This prior has been used previously in
regression settings under the name of smooth-Lasso~\cite{hebiri2011}.
Estimating the model from the data thus consists of minimizing the 
following criterion:\\
\eqcl{
\Jc(\Ub_s, \Vb_s, \Vb) = \sum_{s=1}^S \bigpth{\| \Yb_s - \Ub_s \Vb_s\T \|^2 + \mu \| \Vb_s - \Vb \|^2}\notag}\\
\eqcl{\hspace*{0cm}+ \lambda \, \bigpth{\|\Vb\|_1 + \Vb\T \Lb\, \Vb/2}}
\noindent where, $\|\Vb\|_1$ is the $\ell_1$ norm of $\Vb$, \emph{i.e} the sum
the absolute values, $\Lb$ is the image Laplacian -- $\Vb\T \Lb\, 
\Vb$ is the norm of the gradient. $\lambda$ is a parameter controlling
the amount of prior set on the maps, and thus the amount of sparsity,
that is set by Cross-Validation~(CV). $\mu$ is a parameter controlling the
amount of inter-subject validation, that is set by comparing
intra-subject variance in the observations with inter-subject variance.
For more details about the estimation procedure or the parameter setting,
we refer the reader to~\cite{Varoquaux11}.

\subsection{Resting state MSDL maps}

rs-fMRI runs were analyzed for $S=12$ subjects, consisting of $n=820$ volumes~(time points)
with a $3\,\text{mm}$ isotropic resolution, corresponding to approximately $p=50\,000$ voxels within the
brain. The automatic determination rule of the number of maps exposed in~\cite{Varoquaux10b}
converges to $K=42$. Also, the CV procedure gives us the best CV criterion for $\lambda=2$. 
The group-level maps $\Vb$ are shown in Fig.~\ref{fig:gp_MSDL_maps}. They have been manually classified
in three groups: \emph{Functional}~(F), \emph{Artifactual}~(A) and \emph{Undefined}~(U) maps that appear color-coded
in red, blue and green, respectively.
The undefined class appeared necessary to introduce some confidence measure in our classification and disambiguate
well-established networks~(e.g., dorsal attentional network) from inhomogenous components mixing
artifacts with neuronal regions ~(e.g. like in $\vb_{9}$).
The anatomo-functional description of these group-level maps and their class assignment is given
in Table~\ref{tab:table_MSDL_maps}.
The same rules applied for individual maps $\Vb_s$. In what follows, we will denote by \Fc, \Ac and \Uc
the index sets of F/A/U-maps, respectively and by $\card\pth{\Fc}=25$, $\card\pth{\Ac}=13$ and $\card\pth{\Uc}=4$
their respective size.

\begin{figure*}
\bcc
\begin{tabular}{@{}c@{}c@{}c@{}c@{}c@{}c@{}c}
\encFunc{\figc[width=2.cm]{figure1/sp_map00}}&\encFunc{\figc[width=2.cm]{figure1/sp_map01}}&
\encArt{\figc[width=2.cm]{figure1/sp_map02}}&
\encFunc{\figc[width=2.cm]{figure1/sp_map03}}& \encFunc{\figc[width=2.cm]{figure1/sp_map04}}&
\encFuz{\figc[width=2.cm]{figure1/sp_map05}}&\encArt{\figc[width=2.cm]{figure1/sp_map06}}\\
\encFunc{\figc[width=2.cm]{figure1/sp_map07}}&
\encFuz{\figc[width=2.cm]{figure1/sp_map08}}& \encFunc{\figc[width=2.cm]{figure1/sp_map09}}&
\encArt{\figc[width=2.cm]{figure1/sp_map10}}&\encFunc{\figc[width=2.cm]{figure1/sp_map11}}&
\encFunc{\figc[width=2.cm]{figure1/sp_map12}}&
\encArt{\figc[width=2.cm]{figure1/sp_map13}}\\
\encFunc{\figc[width=2.cm]{figure1/sp_map14}}&
\encFunc{\figc[width=2.cm]{figure1/sp_map15}}&\encArt{\figc[width=2.cm]{figure1/sp_map16}}
&\encFunc{\figc[width=2.cm]{figure1/sp_map17}}&
\encFunc{\figc[width=2.cm]{figure1/sp_map18}}& \encArt{\figc[width=2.cm]{figure1/sp_map19}}&
\encFunc{\figc[width=2.cm]{figure1/sp_map20}}\\
\encFunc{\figc[width=2.cm]{figure1/sp_map21}}&\encFunc{\figc[width=2.cm]{figure1/sp_map22}}&
\encArt{\figc[width=2.cm]{figure1/sp_map23}}& \encFunc{\figc[width=2.cm]{figure1/sp_map24}}&
\encFunc{\figc[width=2.cm]{figure1/sp_map25}}&\encArt{\figc[width=2.cm]{figure1/sp_map26}}&
\encFunc{\figc[width=2.cm]{figure1/sp_map27}}\\
\encFuz{\figc[width=2.cm]{figure1/sp_map28}}& \encArt{\figc[width=2.cm]{figure1/sp_map29}}&
\encFunc{\figc[width=2.cm]{figure1/sp_map30}}&\encFunc{\figc[width=2.cm]{figure1/sp_map31}}&
\encFunc{\figc[width=2.cm]{figure1/sp_map32}}&
\encFunc{\figc[width=2.cm]{figure1/sp_map33}}& \encArt{\figc[width=2.cm]{figure1/sp_map34}}\\
\encArt{\figc[width=2.cm]{figure1/sp_map35}}&\encFunc{\figc[width=2.cm]{figure1/sp_map36}}&
\encFunc{\figc[width=2.cm]{figure1/sp_map37}}&
\encFunc{\figc[width=2.cm]{figure1/sp_map38}}& \encFuz{\figc[width=2.cm]{figure1/sp_map39}}&
\encArt{\figc[width=2.cm]{figure1/sp_map40}}&\encArt{\figc[width=2.cm]{figure1/sp_map41}}
\end{tabular}\vspace*{-.5cm}
\ecc
\caption{From left to right and top to bottom, group-level MSDL maps $\Vb=\cro{\vb_1\I\cdots \I \vb_{42}}$ inferred from the
multisubject~($S=12$) resting-state fMRI dataset~(Neurological convention: left is left).
Functional~(F), Artifactual~(A) and Undefined~(U) maps appear color-coded boxes in \textcolor{red}{red},
\textcolor{blue}{blue} and \textcolor{green}{green}, respectively.  Let us denote \Fc, \Ac and \Uc
the index sets of F/A/U-maps, respectively and $\card\pth{\Fc}=25$, $\card\pth{\Ac}=13$ and $\card\pth{\Uc}=4$
their respective size. Each map $\vb_k$ consists of loading parameters within the $(-1,1)$ range where positive and negative values
are depicted by the hot and cold parts of the colorbar.\label{fig:gp_MSDL_maps}}
\end{figure*}

\begin{table}%
\caption{Classification of group-level $\Vb=\cro{\vb_1\I\cdots\I\vb_{42}}$ maps according to the F/A/U labelling. The F-maps have
been subdivided in different functional networks: Attentional, Default Mode Network, Motor, Visual.
Basal Ganglia~(Thalamus, Caudate and Putamen) and cerebellum have been put together under the Non cortical label.
They will be considered together in the following set: $\Nc=\acc{\text{Att}, \text{DMN}, \text{Mot}, \text{N-c},\text{Vis}}$.
The artifacts have been distinguished in four types: Ventricles, White Matter, Movement and Other. The corresponding
set will be denoted $\Tc =\acc{\text{Ven}, \text{WhM}, \text{Mov}, \text{Oth}}$.
\label{tab:table_MSDL_maps}}
\vspace*{-.25cm}
\bcc
{\small
\begin{tabular}{@{}c@{}p{5cm}c@{}c}\hline\hline\\
{\small\bf  Index} & {\small\bf Anatomo-functional description} & {\small\bf Label} & {\small \bf Network}\\\hline
$\vb_1$ & Ventral primary sensorimotor cortex & F ($\fb_1$)& Mot. \\
$\vb_2$ & Dorsal primary motor cortex or edge of recorded volume & U ($\ub_1$)& \\
$\vb_3$ & Midbrain & A ($\ab_1$)& Oth.\\
$\vb_4$ & Precuneus, posterior cingulate cortex & F ($\fb_2$)& DMN \\
$\vb_5$ & Calcarine cortex (V1) & F ($\fb_3$)& Vis. \\
$\vb_6$ & Anterior cerebellar lobe  & F ($\fb_4$)& N-c \\
$\vb_7$ & Ventricles & A ($\ab_2$)& Ven.\\
$\vb_8$ & Caudate, Thalamus and Putamen &F ($\fb_5$)& N-c\\
$\vb_9$ & Pre- and supplementary motor cortex & U ($\ub_2$) &\\
$\vb_{10}$& Occipital cortex &F ($\fb_6$)& Vis. \\
$\vb_{11}$&  Ventricles &A ($\ab_3$)& Ven.\\
$\vb_{12}$&  Median prefrontal cortex & F ($\fb_7$)& DMN\\
$\vb_{13}$& Right lateralized fronto-parietal cortex & F ($\fb_8$)& Fr.-par. \\
$\vb_{14}$& Ventricles & A ($\ab_4$)& Ven.\\
$\vb_{15}$& Superior temporal and inferior frontal gyrus & F ($\fb_9$)& Lang. \\
$\vb_{16}$& Primary sensorimotor cortex & F ($\fb_{10}$)& Mot. \\
$\vb_{17}$& artifact & A ($\ab_5$)&  Oth.\\
$\vb_{18}$& Dorsal occipital cortex &F ($\fb_{11}$)& Vis.\\
$\vb_{19}$& Supratemporal cortex &F ($\fb_{12}$)&Aud.\\
$\vb_{20}$& Semioval center (white matter) & A ($\ab_6$)& WhM.\\    
$\vb_{21}$& Anterior insula and cingulate cortex & F ($\fb_{13}$)&\\
$\vb_{22}$& Frontal Eye Fields~(FEF), intra-parietal cortex & F ($\fb_{14}$)& Att. \\
$\vb_{23}$& Ventral occipital cortex &F ($\fb_{15}$)&Vis. \\
$\vb_{24}$& Semioval center (white matter) & A ($\ab_7$)& WhM.\\
$\vb_{25}$& Lateral occipital cortex &F ($\fb_{16}$)&Vis. \\
$\vb_{26}$& Parieto-occipital cortex & F ($\fb_{17}$)& Vis. \\
$\vb_{27}$& Extracerebral space & A ($\ab_8$)& Oth.\\
$\vb_{28}$& Left lateralized ventral fronto-parietal cortex & F ($\fb_{18}$)& Fr.-par.\\ 
$\vb_{29}$& Retrosplenial and anterior occipital cortex & U ($\ub_{3}$)&\\
$\vb_{30}$& White matter & A ($\ab_9$)& WhM.\\
$\vb_{31}$& Left lateralized fronto-parietal system & F ($\fb_{19}$)& Fr.-par. \\
$\vb_{32}$& Right lateralized ventral fronto-parietal system & F ($\fb_{20}$)& Att. \\
$\vb_{33}$& Mesial temporal system & F ($\fb_{21}$)&\\ 
$\vb_{34}$& Dorsomedian frontal cortex & F ($\fb_{22}$)& DMN \\
$\vb_{35}$& White matter & A ($\ab_{10}$)& WhM.\\
$\vb_{36}$& Motion-related artifact & A ($\ab_{11}$)& Mov.\\
$\vb_{37}$& Bilateral prefrontal cortex and anterior Caudate &F ($\fb_{23}$)& \\
$\vb_{38}$& Left lateralized temporo-parietal junction and inferior frontal gyrus &F ($\fb_{24}$)& Att.\\
$\vb_{39}$& Right lateralized temporo-parietal junction and inferior frontal gyrus &F ($\fb_{25}$)& Att.\\
$\vb_{40}$& Bilateral superior parietal lobe &U ($\ub_{4}$)\\
$\vb_{41}$& White matter &A ($\ab_{12}$)& WhM.\\
$\vb_{42}$& artifact &A ($\ab_{13}$)& Oth.
\end{tabular}
}
\ecc
\end{table}

To compare spontaneous and evoked activity, the same spatial decomposition was used on resting-state~(run~1, Rest) and
task-related data, which were acquired during an auditory detection task~(run~2, Task).
In practice, this consists of projecting the task-related fMRI data $\wt{\Yb}_s$ onto the inferred spatial maps $\Vb_s$
by minimizing the following least square criterion, $\|\wt{\Yb}_s -\Wb_s \Vb_s\T\|^2$, with respect to $\Wb_s$. The
time series solution admits a closed-form expression: $\wt{\Ub}_s = \wt{\Yb}_s \Vb_s\bigpth{\Vb_s\T\Vb_s}\M$.
The subsequent scale-free analysis is applied to the two sets of $n \times K$ map-level fMRI
time series $\Ub_s=\cro{\ub_{s,1}\I\ldots\I\ub_{s,K}}\T$ and $\wt{\Ub}_s=\cro{\wt{\ub}_{s,1} \I\ldots\I\wt{\ub}_{s,K}}\T$
in a univariate manner, that is to each time series $\ub_{s,k}$ and $\wt{\ub}_{s,k}$ for Rest and Task, respectively.

\section{Scale-free: Intuition, models and analyses}
\label{sec:scaling}

\subsection{Intuition}

In the analysis of evoked brain activity, it is common to seek correlations of BOLD signals with
any a priori shape of the hemodynamic response convolved with the experimental paradigm. In the frequency
domain, this amounts to seeking 
response energy concentration in pre-defined spectral bands, 
as induced for instance by periodic stimulation~(e.g. flashing checkerboards). 
In resting-state fMRI, it is now well admitted that intrinsic brain activity is characterized
by scale-free properties~\cite{Zarahn97,He11}. This constitutes a major change in paradigm as it implies
that brain activity is not to be analyzed via the amounts of energy it shows within specific and a priori chosen
frequency bands, but instead via the fact that all frequencies are jointly contributing in an equivalent manner to
its dynamics. \emph{Scale-free} dynamics are usually described in the spectral domain by a power-law decrease: 
Let $Y(t)$ denote the signal quantifying brain activity and $\Gamma_Y(f)$ its Power Spectral Density~(PSD).
Scale-free property is classically envisaged as: 
\begin{align}\label{eq-PSD}
\Mc_0: \, & \Gamma_Y(f) \simeq C |f|^{-\beta}, \, \beta \geq 0,  
\end{align}
with $ f_m \leq |f| \leq f_M, \, \,  f_M/f_m \gg 1 $.
Such a power law beha\-vior over a broad range of frequencies implies that no frequency in that range plays a specific role,
or equivalently, that they are all equally important. 
To analyze brain activity, this power law relation thus becomes a more important feature than the energy measured
at some specific frequencies.
For instance, it implies that energy at frequency $f_1$ can be deduced from energy at frequency $f_2$ according
to~\cite{He11}: 
\begin{align}\label{eq-PSDb}
\Gamma_Y(f_2) &= \Gamma_Y(f_1)  \left( |f_2|/|f_1| \right)^{-\beta}.
\end{align}
In the scale-free framework, one therefore tries to quantify brain activity by considering the
scaling exponent $\beta$~(or variants) as the key descriptor.
Let us moreover note that the terminology \emph{scale-free} is equivalent to \emph{scale invariance}
or simply \emph{scaling}, encountered in other scientific fields, where this property has also been found to play a
central role~(cf.~\cite{abfrv02,Ciuciu08,ABRY:2012:A}).

\subsection{Scale-free models}\label{sec-scalefreemodel}

\subsubsection{From Spectrum to Increments}
Though appealing, Eqs.~\eqref{eq-PSD}--\eqref{eq-PSDb} do not provide practitioners with a versatile enough
definition of \emph{scale-free} with respect to real-world data analysis.
Indeed, they concentrate only on the second
order statistics and hence account neither for the marginal distribution~(first order statistics) of the signal $Y$,
nor for its higher order dynamics (or dependence structure).
For instance, it does not indicate whether data are jointly Gaussian or depart, weakly or
strongly, from Gaussianity. 

To investigate how to enrich Model $\Mc_0$, let us assume for now that $Y$ consists of a stationary jointly
Gaussian process, with PSD as in Eq.~\eqref{eq-PSD}.
Equivalently, this implies that the covariance function behaves as  $C_Y(\tau) \sim \sigma^2_Y (1+C' |\tau|^{-\alpha}) $, for $  \tau_m \leq \tau \leq \tau_M$, with $\alpha = 1- \beta$. 
A simple calculation hence shows that $\E  (Y(t+\tau) -Y(t))^2 = \E  Y(t+\tau)^2  + \E Y(t)^2  -2 \E Y(t+\tau)Y(t) =  c_2 |\tau|^{-\alpha}$. 
The Gaussianity of $Y$ further implies that $\,\forall q > -1$: 
\begin{align}\label{eq-incY}
\E |Y(t+\tau) -Y(t)|^q &= c_q |\tau|^{-\frac{q\beta}{2}},� \, \,  \tau_m \leq \tau \leq \tau_M. 
\end{align}

Defining $X(t) = \int^t Y(s) ds$, Eq.~\eqref{eq-incY} straightforwardly implies that, as long as
$\tau_m \leq \tau_1, \tau_2 \leq \tau_M$:
\begin{align}
\label{eq-fdd}
\bigacc{\frac{X(t+\tau_1)-X(t)}{\tau^H_1}}_{t \in \R} \stackrel{fdd}{=}
\bigacc{\frac{X(t+\tau_2)-X(t)}{\tau^H_2}}_{t \in \R}, 
\end{align}
with $H = (-\alpha/2) = (\beta+1)/2$, and where $ \stackrel{fdd}{=}$ means equality of all joint finite dimensional
distributions: i.e.,  $(X(t+\tau_1)-X(t))/\tau^H_1$ and
$(X(t+\tau_2)-X(t))/\tau^H_2$ have the same joint distributions. 
In turn, this implies that $\,\forall q > -1$,
such that $\E |X(t)|^q < \infty$: 
\begin{align}
\label{eq-inc}
\E |X(t+\tau) -X(t)|^q &= c_q |\tau|^{qH}, \, \,  \tau_m \leq \tau \leq \tau_M, \makebox{ or }
\end{align}
\begin{equation}
\label{eq-incb}
\E |X(t+\tau_2) -X(t)|^q = \E |X(t+\tau_1) -X(t)|^q  \left(\frac{|\tau_2|}{|\tau_1|}\right)^{qH},
\end{equation}
with $ \tau_m \leq \tau_1, \tau_2 \leq \tau_M$, which are reminiscent of Eqs.~\eqref{eq-PSD}--\eqref{eq-PSDb}.

\subsubsection{Self-Similar processes with stationary increments}
Eqs.~\eqref{eq-inc}--\eqref{eq-incb} turn out to hold not only for jointly Gaussian $1/f$-processes but for a much wider
and better defined class, that of self-similar processes with stationary increments, referred to as $H$-sssi processes,
and defined as~(cf.~\cite{Samorodnitsky1994}): 
\begin{align}
\label{eq-Hsssi}
\Mc_1: &\quad \{X(t)\}_{t \in \R} \stackrel{fdd}{=} \{a^H X(t/a)\}_{t \in \R},
\end{align}
 $\forall a > 0$,  $H \in (0,1)$. 
Essentially, it means that $X$ cannot be distinguished (statistically) from any copy, dilated by
scale factor $a>0$, on condition that the amplitude axis is scaled by $a^H$. 
Parameter $H$ is referred to as the self-similarity exponent.
A major practical consequence of this definition consists of the fact that Eqs.~\eqref{eq-inc}--\eqref{eq-incb} hold for all
$\tau $ (resp., $\tau_1, \tau_2$).

The central benefit of such a definition is that it does not require the data to be Gaussian but provides both
theoreticians and practitioners with a well-defined model.
For analysis, fMRI data can hence be envisaged as the increment process $Y(t) = X(t+\tau_0)-X(t)$ of an $H$-sssi
process $X$ (where $\tau_0$ is an arbitrary constant chosen to make sense with respect to physiology and data
acquisition set up, e.g. $\tau_0=\text{TR}$).
This constitutes a second model to account for scale-free properties in data, that encompasses the \emph{simpler}
 $1/f$-spectrum first model.

Further, if joint Gaussianity is assumed, the model becomes even more precise as the only Gaussian $H$-sssi process $X$ is
the so-called fractional Brownian motion (fBm), cf. e.g.,~\cite{mvn68}, hereafter labelled $X(t) \equiv B_H(t)$. 
The corresponding increment process $Y(t) = G_H(t) = B_H(t+1)-B_H(t)$ is termed fractional Gaussian noise (fGn). 
Additionally, note that it may sometimes constitute a practical and relevant challenging issue to decide whether brain
activity is better modelled by the  $H$-sssi process $X$ (hence a non stationary process) or by its
increment process $Y$ (hence a stationary process) (cf. e.g.,~\cite{Ciuciu08,He10,He11}).

\subsubsection{Multifractal processes}

In a number of situations, it has been actually observed on a variety of real-world data of very different
nature~(cf. e.g., \cite{abfrv02,ABRY:2012:A} for reviews) that Eq.~\eqref{eq-inc} holds over a wide range of $\tau$s,
however, with scaling exponents that depart significantly from the theoretical linear behavior $qH$:
\begin{align}
\label{eq-incc}
\E |X(t+\tau) -X(t)|^q &= c_q |\tau|^{\zeta(q)}, \, \,  \tau_m \leq \tau \leq \tau_M.
\end{align}
The generic behaviors modeled by Eq.~\eqref{eq-incc} can be considered as a practical or operational,
definition of scale-free property. 
Let us note that, by nature, $\zeta(q)$ is necessarily a concave function of $q$ (cf. e.g., \cite{WENDT:2007:E}).

Scaling exponents $\zeta(q) $ that are strictly concave rule out the use of $H$-sssi process as models. 
Instead, a broader class should be used, referred to as that of \emph{multifractal processes}.
This is however a large and not-well defined class of processes. 
For the purposes of this contribution, let us use a particular subclass of multifractal processes defined as fBm
subordinated to a \emph{multiplicative  Compound Poisson cascade}:
\begin{align}
\label{eq-MF}
\Mc_2:&\quad X(t) : = B_H(A(t)), \WHERE A(t) = \int^t W(s) ds, 
\end{align}
with $W(s) $ a multiplicative Compound Poisson cascade (or martingale), such as those defined in \cite{BaM02}. 
The complete definition of these cascades has been given and studied with details elsewhere and is hence not recalled
here (cf. \cite{BaM02,Bacry2001,Chainais2005}). 
It is enough to emphasize that they rely on the choice of positive random variables whose
moments of order $q$ define the $\zeta(q)$. 
The process $X$ thus defined satisfies Eq.~\eqref{eq-incc} with strictly convex tunable scaling exponents $\zeta(q)$,
has stationary increments $Y$, and has distributions that depart from strict jointly Gaussian laws. 
Such departures, that may however turn subtle and hard to detect in practice, are precisely quantified by the departure of
$\zeta(q)$ from a linear behavior in $q$.
The $\zeta(q)$ therefore convey a rich information about data $X$, and hence about $Y$, as they account for the entire
dependence structure of the data, hence both to the time dynamic and distributions of data.
Their accurate estimation from real-world data therefore naturally constitutes an important practical challenge discussed
below.

\subsection{Scale-free analysis}

\subsubsection{From spectrum to wavelet analysis}
Assuming that data $Y$ have a power-law spectrum beha\-vior as in Eq.~\eqref{eq-PSD}, it is natural to rely on spectral
estimation to measure $\beta$. 
A classical tool in spectrum analysis is the Welch estimator that consists in splitting data $Y$ into blocks
and in averaging the squared Fourier transforms computed independently over each block. 
For scale free data, it is hence expected that: 
\begin{align}\label{eq-PSDestim}
\hat \Gamma_Y(f) &= \sum_k |\langle Y, g_{f,k} \rangle|^2 \simeq C |f|^{-\beta},
\end{align}
where the $g_{f,k} = g_0(t-k) e^{\imath 2 \pi ft}$ are translated into time and into frequency templates of a reference
pattern $g_0(t)$. 
This relation can be further used to estimate $\beta$. 

It has been shown that wavelet transforms can achieve better performance both in the analysis of scale-free properties
in real-world data, and in the estimation of the corresponding scaling parameters (cf.~\cite{agf95,avf98,va01}). 
The discrete wavelet transform~(DWT) coefficients of $Y$ are defined as: 
\begin{align}\label{eq-wc}
d_Y(j,k) &= \displaystyle\int_{\R} Y(t)  \;  2^{-j}\psi_0(2^{-j}t-k) \, dt  \equiv \langle Y, \psi_{j,k} \rangle, 
\end{align}
where the $\psi_{j,k} = 2^{-j}\psi_0(2^{-j}t-k)$ consists of templates of a reference pattern $\psi_0$
translated in time and dilated~(by a factor $a=2^j$). It is referred to as the {\em mother-wavelet}:
an elementary function, characterized by fast exponential decays in both the time and frequency domains, as well as 
by a strictly positive integer $N_\psi \geq 1$, the {\em number of vanishing moments},
defined as $ \forall k = 0,1,\ldots, N_\psi-1$, $ \int_\R t^k\psi_0(t) dt \equiv 0$ and $ \int_\R t^N \psi_0(t) dt \neq 0$. 
Note the choice of the $L^1$-norm (as opposed to the more common $L^2$-norm choice) that better matches sca\-ling analysis. 
For further introduction to wavelet transforms, the reader is referred to e.g.,~\cite{Mallat2009}.

Defining $S^d_Y(j,2) = { 1 \over n_j } \sum_{k=1}^{n_j} |d_Y(j,k)|^2$ (with $n_j$ the number of $d_X(j,k)$
available at scale $2^j$),  one obtains (cf. \cite{agf95}): 
\begin{align}
\label{eq-wavsp}
\E S^d_Y(j,2) &=  \displaystyle\int_{\R}  \Gamma_Y(f) |{\Psi}_0(2^j f)|^2 df
\end{align}
where ${\Psi}_0$ denotes the Fourier transform of $ \psi_0$. 
This indicates that $S^d_Y(j,2) $ can be read as a wavelet based estimate of the PSD and is hence
referred to as the \emph{wavelet spectrum}. It measures the amount of energy of $Y$ around the frequency $f_j=f_0/2^j$
where $f_0$ is a constant that depends on the explicit choice of $\psi_0$~(for the Daubechies wavelet used here,
$f_0\simeq 3 f_s /4$
with $f_s$ the sampling frequency).
This correspondence between the Fourier and wavelet spectra is illustrated on fMRI signals in Fig.~\ref{fig:spectra}. 
For scale-free processes satisfying Eq.~\eqref{eq-PSD}, it implies:
\begin{align*} 
S^d_Y(j,2) \equiv { 1 \over n_j } \sum_{k=1}^{n_j} |\langle Y, \psi_{j,k} \rangle|^2 
\simeq C_2 2^{j (\beta-1)}, \, \, a_m \leq 2^j \leq a_M. 
\end{align*}
While this formally looks like Eq.~\eqref{eq-PSDestim}, it has been shown
in detail how and why the wavelet spectrum yields better estimates of the scaling exponents $\beta$ than Welch based-ones,
both in terms of estimation performance and robustness to various forms of non-stationarity
in data that may be confused with scale-free behaviors~\cite{agf95,avf98,va01}. 
Notably, it was shown how wavelet analysis enables to disentangle non stationarity, stemming from fMRI environment,
from true long memory in brain activity.
Also, the wavelet spectrum avoids the potentially difficult issue that consists of deciding a priori whether
empirical data are better modeled by $Y$ or $X$, needed by classical spectrum estimation, that can only be applied to
stationary data. 
In a nutshell, these benefits stem from the use of the change of scale operator to design the analysis tool,
that intuitively matches scale-free behavior more naturally than a frequency shift operator.

\subsubsection{From 2nd to other statistical orders: Wavelet leaders}

As discussed in Section~\ref{sec-scalefreemodel}, analyzing in-depth scale free pro\-perties implies investigating
not only the spectrum~(i.e., the se\-cond order statistics of data) but rather the entire dependence structure,
i.e., the whole range of avai\-lable statistical
orders $q$. It had initially been thought that this would amount to extending the definition of
$S^d_Y(j,2)$ to other orders $q$, $  S^d_Y(j,q) \equiv { 1 \over n_j} \sum_{k=1}^{n_j} |\langle Y, \psi_{j,k} \rangle|^q$. 
It has howe\-ver recently been shown that this approach, though intuitive and appealingly \emph{simple},
fails to yield satisfactory estimation of the $\zeta(q)$.
Notably, wavelet coefficients show little power in enabling practitioners to decide whether $\zeta(q)$ is a linear
or strictly concave function
of $q$. Instead, it is now well documented that the estimation of the $\zeta(q)$ should be based on Wavelet
Leaders~\cite{WENDT:2007:E}.

Let us now assume that $\psi_0$ has a compact time support and introduce the global regularity of $Y$, $h_m$, defined as:
$h_m = \liminf_{2^j \rightarrow 0}\froc{\log\bigpth{\sup_{k} |d_Y(j,k)|}}{\log ( 2^j ) }$.
Therefore, $h_m $ can be estimated by a linear regression of the log of
the magnitude of the largest wavelet coefficient at scales $2^j$ versus the log of the scales
$ 2^j$~\cite{WENDT:2007:E,ABRY:2012:A}.
Let $\gamma \geq 0$ be defined as, with $ \epsilon > 0$: $ \gamma = 0$ if $ h_m > 0$, and $\gamma = - h_m + \epsilon$
otherwise. 
Further, let $\lambda_{j,k} $ denote the dyadic interval
$\lambda_{j,k}=[k 2^j,(k +1)2^j)$, and denote by $3 \lambda_{j,k}$ the union of
$\lambda_{j,k}$ and its 2 closest neighbours,
$3\lambda_{j,k}= [(k-1) 2^j,(k+2)2^j)$.
The wavelet leaders $L^{(\gamma)}_Y $ are defined as
$L^{(\gamma)}_Y(j,k)=\sup_{\lambda'\subset3\lambda_{j,k}} 2^{\gamma j} |d_Y(\lambda')| $. 
In practice, $L^{(\gamma)}_Y(j,k)$ simply consists of any of the largest coefficients
$2^{\gamma j} |d_Y(\lambda')|$ located at scales finer or equal to $2^j$ and within a small time neighborhood.
It is then necessary to form the so-called wavelet Leader structure functions that reproduce the scale-free
properties in $Y$ according to:
\begin{align}\label{eq-leaderscaling}
 S^L_Y(j,q,\gamma) &\equiv { 1 \over n_j } \sum_{k=1}^{n_j} (L^{(\gamma)}_Y(j,k))^q 
 \simeq c_q 2^{j\zeta(q, \gamma)},
\end{align}

\noindent Moreover, for a large class of processes, one has:
 $\zeta(q, \gamma) = \zeta(q) +\gamma  q $. 
For all real-world data analyzed so far with WLMF, this relation is found to hold, by varying
$\gamma$ (cf.~\cite{WENDT:2007:E,ABRY:2012:A} for a thorough discussion).
This has also been verified empirically for fMRI data.
Further, because it can take any concave shape, the function $\zeta(q,\gamma)$ is often written as a
polynomial expansion~\cite{aadmv02}: 
$\zeta(q,\gamma) = \sum_{p \geq 1} c^{(\gamma)}_p \froc{q^p}{ p~!}$.
Notably, the second order truncation $\zeta(q, \gamma) \simeq c^{(\gamma)}_1 q + c^{(\gamma)}_2 q^2/2 $
(with $c^{(\gamma)}_2 \leq 0$ by concavity) can be regarded as a potentially interesting approximation that captures
the crucial information regarding whether the $\zeta(q, \gamma)$ are linear in $q$ (hence indicating $H$-sssi models)
or strictly concave (hence suggesting multiplicative cascade models).  
Interestingly, the coefficients $ c^{(\gamma)}_{p} $ entering the polynomial expansion of
$\zeta(q,\gamma)$ are not abstract figures but rather turn out to be quantities deeply tied to
the scale-free properties of $Y$, as they are related to the scale dependence of the cumulants of order
$p\geq 1$, $C^{(\gamma)}_Y(j,p) $,  of the random variable $ \ln L^{(\gamma)}_Y(j,k)$: 
\bal
\label{kappa_form}
\forall p \geqslant 1, C^{(\gamma)}(j,p)_Y &= c^{(\gamma)}_{0,p}+c^{(\gamma)}_p\ln 2^j.
\eal
Eqs.~\eqref{eq-leaderscaling}--\eqref{kappa_form} suggest that the $\zeta(q,\gamma) $ or $ c^{(\gamma)}_p $
can be efficiently estimated from linear regressions: $\hat \zeta(q,\gamma)  =  \sum_{j=j_1}^{j_2}w_j \log_2 S^L_Y(j,q,\gamma)$ and 
$\hat c^{(\gamma)}_p = \log_2e \sum_{j=j_1}^{j_2}w_j\hat C^L_Y(j,p,\gamma)$. 
The weights $w_j$ are chosen to perform ordinary~(or non weighted) least squares estimation~(cf.~\cite{va01} for discussion). 
Further,  $\zeta(q, \gamma) = \zeta(q) +\gamma  q $ obviously implies that $ c_{1} = 
c^{(\gamma)}_1 - \gamma$  and $\forall p \geq 2$,  $c_{p} = c^{(\gamma)}_p$. 

This wavelet-Leader based analysis of scale-free properties is intimately and ultimately related to
\emph{multifractal analysis}, the detailed introduction of which is beyond the scope of the present contribution.
We restate here only its essence.
Multifractal analyses describe globally the fluctuations along time of the local regularity of a signal $Y(t)$. 
This local regularity is measured by the so-called H\"older exponent $h(t) $, that essentially compares $Y$ around time $t_0$
against a local power-law behavior: 
$ |Y(t)-Y(t_0)| \leq |t-t_0|^h$, $  |t-t_0| \rightarrow 0$. 
The variations of $h$ along time are then described globally via the multifractal spectrum, consisting of the collection of
Hausdorff dimensions, ${\cal D}(h)$,  of the sets of points $\{t, h(t) = h\}$. 
In practice, the multifractal spectrum is estimated indirectly via (a Legendre transform of) the function $\zeta(q)$. 
The approximation  $\zeta(q) \simeq c_1 q + c_2 q^2/2 $ translates into ${\cal D}(h) \simeq 1 - (h-c_1)^2/(2|c_2|)$.
For thorough and detailed introductions to multifractal analysis, the reader is referred to e.g., \cite{WENDT:2007:E}. 
Examples of such multifractal spectra estimated using the WLMF from real fMRI signals are illustrated
in Fig.~\ref{fig:spectra}(b).
An outcome of the mathematical theory underlying multifractal analysis, of key practical importance and impact,
is the following: 
the function ${\cal D}(h)$ theoretically constitutes a rich characterization of the scale-free properties
of a signal $Y$ and its complete and entire estimation requires the use, in Eq.~\eqref{eq-leaderscaling},
of both positive and negative order $q$s, concentrated left and right around $0$ \cite{WENDT:2007:E}. 

\section{Multifractal analysis of MSDL maps}\label{sec:results}

\subsection{Single subject analysis}

\subsubsection{Scaling range}

For analysis, orthonormal minimal-length time support Daubechies's wavelets were used with $N_\psi =3$. 
Scale-free properties are systematically found to hold within a 4-octave range~($(j_1,j_2)=(3,6)$), corresponding to a
frequency range of $[0.008, 0.063] $ Hz\footnote{The scale and band-specific central frequency are
related according to $f_j=3 f_e/(4 2^j)$.}, which is hence consistent with the upper limit $0.1$Hz classically
associated with the hemodynamics boundary and scaling in fMRI data~\cite{Cordes01}.

\subsubsection{Fourier vs. wavelet spectra}

\begin{figure}
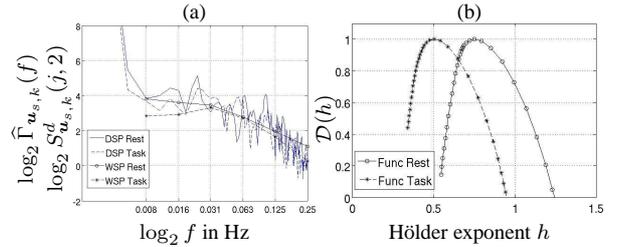

\bcc
\btabu{@{}c@{}c@{}c@{}c}
\yaxis{\btabu{c}{\footnotesize $\log_2 \wh{\Gamma}_{\ub_{s,k}}(f)$} \\
{\footnotesize $\log_2 S^d_{\ub_{s,k}}(j,2)$}\etabu}&
\figc[height=2.5cm]{figure2/DSP_28_cropped} &
\yaxis{{\footnotesize $\Dc(h)$}}&
\figc[height=2.5cm]{figure2/Dh_28_cropped}\\[-3.cm]
& {\footnotesize (a)} & &{\footnotesize (b)}\\[2.5cm]
& {\footnotesize $\log_2 f$ in Hz} & 
& {\footnotesize H\"older exponent $h$ }\\
\etabu\vspace*{-.6cm}
\ecc
\caption{{{\bf (a):} Welch~(\textcolor{blue}{blue} curves) vs. Wavelet~({\bf \textcolor{black}{black}} curves)
spectra associated with a F-map~($\fb_{18}$). Solid and dashed lines correspond to rest and task, respectively.
{\bf (b):} Corresponding multifractal spectra $\Dc(h)$. \label{fig:spectra}}}
\end{figure}

For illustrative purposes, two time series corresponding to a functional map~($k\!=\!28$, $\fb_{18}$ in
Tab.~\ref{tab:table_MSDL_maps}), were selected in the rest and task runs from the first subject.
In Fig.~\ref{fig:spectra}(a), the Fourier spectrum
estimate~($\log_2 \wh{\Gamma}_{\ub_{s,k}}(f)$) based on Welch's
averaged periodogram and its wavelet spectrum counterpart ($\log_2 S^d_{\ub_{s,k}}(j,2)$)  are found to closely match,
as predicted by Eq.~\eqref{eq-wavsp}.
Interestingly, Fig.~\ref{fig:spectra}(a) shows that the $\beta$ exponent, measured within frequency range
$[0.008, 0.063] $ Hz,  in Eq.~\eqref{eq-PSD}~(i.e. the neg-slope of
the log-spectra $\log_2 \wh{\Gamma}_{\ub_{s,k}}(f)$)
decreases with task-related activity in $\fb_{18}$. This amounts to observing lower
Hurst exponent $H=(\beta-1)/2$ in the task-related dataset:
$\wh{H}^{\rm R}_{\fb_{18}}\simeq  0.66 $ and $\wh{H}^{\rm T}_{\fb_{18}}\simeq 0.5$.
As shown in the following, this decrease of self-similarity
is not specific to functional maps and will be observed in artifactual and undefined maps.
Following~\cite{He11}, the stationarity of fMRI signals is confirmed since we systematically observed $\wh{H}_k^{\rm R,T}<1$.

\subsubsection{Multifractal spectrum}

For the same time series, MF spectra $\Dc(h)$, estimated using the WLMF tool described above,
are depicted in Fig.~\ref{fig:spectra}(b). 
The decrease of self-similarity between rest and task is captured by a shift to the left
of the position $\wh{c}_1$ of the maximum of ${\cal D}(h)$: $\bigpth{(\wh{c}_1)_{\fb_{18}}^{\rm R},
(\wh{c}_1)_{\fb_{18}}^{\rm T}}\!=\!(0.75,0.5)$
It should also be noted that parameter $\wh{c}_1$ systematically takes values that are close to those of the Hurst exponent. 
This is consistent with the theoretical mode\-ling of scale-free property that establishes a clear connection between $c_1$ and $H$
and predicts $c_1 \!\simeq\! H$ (cf.~\cite{WENDT:2007:E}).
Therefore, in the following, $c_1$ will be referred to as the self-similarity parameter although this is a slight
misnomer.
Further, Fig.~\ref{fig:spectra}(b) confirms the presence of multifractality in fMRI data as strictly negative
$c_2\!<\!0$ are almost always observed. 
Indeed, parameter $c_2$ quantifies the width of ${\cal D}(h)$ (as a curvature radius of~${\cal D}(h)$ around
$\wh{c}_1$): $\wh{c}_2 \!<\! 0$.
Multifractality is however not specific to a given brain state since we measured
$\bigpth{(\wh{c}_2)_{\fb_{18}}^{\rm R}, (\wh{c}_2)_{\fb_{18}}^{\rm T}}\!=\!(-0.07,-0.06)$.
In this example, multifractality, as measured by the width of the multifractal spectra, is decreased
from rest to task. However, opposite fluctuations will be also observed amongst F-maps.

The sole two self-similarity and multifractality parameters $c_1$ and $c_2$ are therefore used from now on as
sufficient and relevant descriptors of the scale-free properties of fMRI signals (superscript $\gamma$ is dropped
for the sake of conciseness, while $\gamma$ has been systematically set to $\gamma =  2$).

\subsection{Group-level analysis}

\subsubsection{Group level scale-free properties}\label{subsec:mean_effects}

Let $c_{i,k}^{j,s}$ denote the $\wh{c}_1$ and $\wh{c}_2$ estimates~(index $i=1\!\!:\!\!2$)
for differents maps~(index $k=1\!\!:\!\!K$), runs~(index $j={\rm R,T}$ for Rest and Task, respectively) and for different
subjects~(index $s$).
The map-dependent group-level values have been computed as $\mu_{i,k}^j\!=\!\sum_{s=1}^S\wh{c}^{j,s}_{i,k}/S$ and
sorted according to their labelling~(F/A/U maps) given in Tab.~\ref{tab:table_MSDL_maps}.
Then, global spatial averaging of the means $\mu_{i,k}^j$ has been performed so as to derive global F/A/U-average
parameter estimates: 
\small$\bar{\mu}_{i,F}^j\!=\!\sum_{k\subset \Fc}\mu_{i,k}^j/\card\pth{\Fc}\, $\normalsize,  $\bar{\mu}_{i,A}^j$ and
$\bar{\mu}_{i,U}^j$ are defined equivalently.
In the same spirit, group-level multifractal attributes $\bar{\mu}_{i,v_\ell}^j$ are derived for each functional
network $v_\ell\in\Nc\!\!=\!\!\acc{\text{Att}, \text{DMN}, \text{Mot}, \text{N-c},\text{Vis}}$ such that
$\bar{\mu}_{i,\nb_\ell}^j\!\!=\!\!\sum_{k\in \nb_\ell} \mu_{i,k}^j/\card\pth{\nb_\ell}$, $\forall \ell=1:5$ and
$j=({\rm R,T})$. We proceed in the same way for analyzing artifact types
$\tb_r\in \Tc\!\!=\!\!\acc{\text{Ven}, \text{WhM}, \text{Mov}, \text{Oth}}$, and computing
$\bar{\mu}_{i,\tb_r}^j$ for $r=1:4$.

\begin{figure}
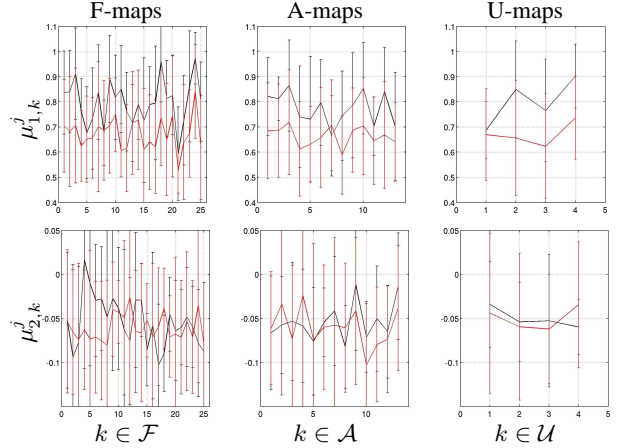

\bcc
\begin{tabular}{@{}c@{}ccc}
\yaxis{$\mu_{1,k}^j$}& \figc[height=2.5cm]{figure3/spDict_c1_N3_j13_j26_Fmaps3} &
\figc[height=2.5cm]{figure3/spDict_c1_N3_j13_j26_Amaps3}&
\figc[height=2.5cm]{figure3/spDict_c1_N3_j13_j26_Umaps3} \\[-3.cm]
 & {\small F-maps} & {\small A-maps}& {\small U-maps}\\[2.65cm]
 \yaxis{$\mu_{2,k}^j$}& \figc[height=2.5cm]{figure3/spDict_c2_N3_j13_j26_Fmaps3} &
 \figc[height=2.5cm]{figure3/spDict_c2_N3_j13_j26_Amaps3}&
 \figc[height=2.5cm]{figure3/spDict_c2_N3_j13_j26_Umaps3}\\[-.1cm]
& {\small $k\in\Fc$} &  {\small $k\in\Ac$} & {\small $k\in\Uc$}
\end{tabular}\vspace*{-.75cm}
\ecc
\caption{From left to right: Group-averaged map-dependent MF parameters $\mu_{1,k}^j$~(top), $\mu_{2,k}^j$~(bottom)
specific to F/A/U-maps defined in Tab.~\ref{tab:table_MSDL_maps}.
\textbf{Black} and \textcolor{red}{red} curves code for $j={\rm R}$~(Rest) and $j={\rm T}$~(Task).\label{fig:errorbarMF_maps}}
\end{figure}

As shown in Fig.~\ref{fig:errorbarMF_maps}[top], the group-averaged values of self-similarity $\mu_{1,k}^j$ lie approximately
in the same range [.55, 1], indicating long memory, for all components~(F/A/U-maps). 
An almost systematic decrease of
self-similarity is observed in the task-related dataset~($\delta_{1,k}\!=\!\mu_{1,k}^{\rm T}-\mu_{1,k}^{\rm R}\!<\!0$), for
$k\in\Fc\cup\Ac\cup\Uc$. This trend is therefore not specific to F-maps.
Moreover, the average decrease computed over F-maps is about the same as the one estimated for
A and U-maps~($\bar{\delta}_{1,F}\!=\!-0.125$, $\bar{\delta}_{1,A}\!=\!-0.11$ and $\bar{\delta}_{1,U}\!=\!-0.13$).
Also, the averaged standard deviations~($\bar{\sigma}^{\rm R}_{1,F}$, $\bar{\sigma}^{\rm R}_{1,A}$ and
$\bar{\sigma}^{\rm R}_{1,U}$) computed over the F/A/U-maps,
are close to each other~$(\bar{\sigma}^{\rm R}_{1,F/A/U}\approx 0.18)$ and systematically increase with
the task-related activity~($\bar{\sigma}^{\rm T}_{1,F/A/U} \!>\! \bar{\sigma}^{\rm R}_{1,F/A/U}$).

Fig.~\ref{fig:errorbarMF_maps}[bottom] illustrates that the group-averaged values of $\mu_{2,k}^{\rm R,T}$ are almost all 
negative in the F/A/U-maps indicating multifractality in fMRI time series irrespective of the map type or
brain state. 
Between rest to task-related situation minor changes
in the A and U-maps are also observed since
$\bars{\delta_{2,k}}\!<\!0.03$ for $k\in\Uc\cup\Ac$  while we measured $\bars{\delta_{2,k}}\!<\!0.08$ for
$k\in\Fc$ ($\delta_{2,k}\!=\!\mu_{2,k}^{\rm T}-\mu_{2,k}^{\rm R}$).
Hence, the level of multifractality does not change much between rest and task in irrelevant maps.
In contrast, large changes in the multifractal parameters are observed in F-maps, while not systematically in
the same direction.
For instance, in cerebellum~($\fb_4$), basal ganglia~($\fb_5$), DMN~($\fb_7$) and fronto-parietal
network~($\fb_8$) evoked activity induces a large 
increase of multifractality~($\delta_{2,k}\!<\!0$) while in the auditory and attentional systems~(e.g. $\fb_{12}$ and
$\fb_{24}$, respectively), which are supposed to be involved in the auditory detection task,
the converse observation holds, i.e. $\delta_{2,k}\!>\!0$. 
Also, it is worth noticing that the averaged standard deviations computed over the A/U-maps increase when switching
 from rest to task~($\bar{\sigma}^{\rm R}_{2,A}\!=\!0.06<\bar{\sigma}^{\rm T}_{2,A}\!=\!0.09$ and
 $\bar{\sigma}^{\rm R}_{2,U}\!=\! 0.05\!<\! \bar{\sigma}^{\rm T}_{2,U}\!=\!0.085$) while they remain at the same
 level in the F-maps: $\bar{\sigma}^{\rm R}_{2,F}\!\approx\! \bar{\sigma}^{\rm T}_{2,F}\!\approx\! 0.08$.

\begin{figure}
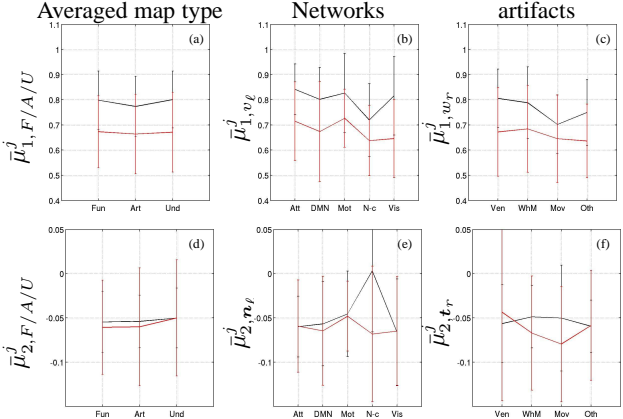

\bcc
\begin{tabular}{@{}c@{}c@{}c@{}c@{}c@{}c}
\yaxis{$\bar{\mu}_{1,F/A/U}^j$}& \figc[height=2.5cm]{figure4/spDict_c1_N3_j13_j26_FAU3} &
\yaxis{$\bar{\mu}_{1,v_\ell}^j$}&\figc[height=2.5cm]{figure4/spDict_c1_N3_j13_j26_Nets3}&
\yaxis{$\bar{\mu}_{1,w_r}^j$}&\figc[height=2.5cm]{figure4/spDict_c1_N3_j13_j26_Arts3} \\[-2.65cm]
& \hspace*{1.75cm}{\tiny (a)} & & \hspace*{1.75cm}{\tiny (b)} & & \hspace*{1.75cm}{\tiny (c)}\\[-.75cm]
&{\small Averaged map type} & & {\small Networks} & & {\small artifacts}\\[2.65cm]
\yaxis{$\bar{\mu}_{2,F/A/U}^j$} & \figc[height=2.5cm]{figure4/spDict_c2_N3_j13_j26_FAU3} &
\yaxis{$\bar{\mu}_{2,\nb_\ell}^j$}& \figc[height=2.5cm]{figure4/spDict_c2_N3_j13_j26_Nets3}&
\yaxis{$\bar{\mu}_{2,\tb_r}^j$}   & \figc[height=2.5cm]{figure4/spDict_c2_N3_j13_j26_Arts3} \\[-2.65cm]
& \hspace*{1.75cm}{\tiny (d)} & & \hspace*{1.75cm}{\tiny (e)} & & \hspace*{1.75cm}{\tiny (f)}\\[2.25cm]
\end{tabular}\vspace*{-.75cm}
\ecc
\caption{
    Group-level MF parameters averaged over the F/A/U-maps $\bar{\mu}_{i,F/A/U}^j$~(Left), the 
    functional networks $\bar{\mu}_{i,\nb_\ell}^j$~(Middle) and the artifact types $\bar{\mu}_{i,\tb_r}^j$~(Right).
    \textbf{Black} and \textcolor{red}{red} curves code for $j={\rm R}$~(Rest) and $j={\rm T}$~(Task).\label{fig:errorbarMF_global}}
    \end{figure}

We computed the grand means of the self-similarity parameters $\bar{\mu}_{1,F/A/U}^{\rm R,T}$
over the F/A/U-maps, respectively, and draw the same conclusion at this macroscopic level, as demonstrated
in Fig.~\ref{fig:errorbarMF_global}(a)-(c): the decrease of self-similarity from rest to task
is not specific to functional components and only slightly fluctuates between networks and artifact types.
Moreover, we did not observe any significant modification of the grand means of multifractal parameter estimates
$\bar{\mu}_{2,F/A/U}^{\rm R,T}$ between rest and task, as illustrated in Fig.~\ref{fig:errorbarMF_global}(d).
This motivated deeper investigations at the network and artifact levels, especially concerning the fluctuation
of multifractality induced by task. Fig.~\ref{fig:errorbarMF_global}(e) reveals that a major increase
of multifractality~($\bar{\mu}_{2,\nb_4}^{\rm T}\! <\! \bar{\mu}_{2,\nb_4}^{\rm R}$) occurred only in the
non-cortical regions while no major change appeared in the artifacts~($\bar{\mu}_{2,\tb_r}^{\rm T}\!\simeq\!
\bar{\mu}_{2,\tb_r}^{\rm R}, \forall r\!\in\!\Tc$) as shown in Fig.~\ref{fig:errorbarMF_global}(f).

\subsubsection{One-sample statistical tests}

To assess the statistical significance of the multifractal parameters for the rest and task-related datasets at
the group-level, we used \emph{one-sided} tests associated with the following null hypotheses
$\forall k\in\Fc\cup\Ac\cup\Uc$:
\begin{equation}
\left.
\begin{array}{ll}
H_{0,j}^{(1,k)}: \mu_{1,k}^{j}&\leqslant 0.5, \quad  \quad\text{(White noise or SRD)} \\
H_{0,j}^{(2,k)}:\mu_{2,k}^{j}&= 0., \quad  \quad\text{(H-sssi process)}.
\end{array}
\right\}
\label{eq:StatsAssumptions}
\end{equation}
We also conducted similar tests at the macroscopic level~($k\in\Nc\cup\Tc$) by replacing $\mu_{i,k}^j$ with
$\bar{\mu}_{i,k}^j$ in the null hypotheses~\eqref{eq:StatsAssumptions}.
Because there is no definite proof nor evidence that MF parameter estimates $\wh{c}_{i,k}^{j,s}$ should be normally
distributed across subjects, we investigated different statistics~(Student-$t$, Wilcoxon's signed rank~(WSR) statistic).
Indeed, other statistics may provide more sensitive results in presence of outliers.
To account for multiple comparisons~($K$ tests performed simultaneously) and to ensure
correct specificity control~(control of false positives), the Bonferroni correction was applied.

Rejecting $H_{0,j}^{(1,k)}$ clearly amounts to localizing brain areas or components eliciting
significant long memory or self-similarity.
Rejecting $H_{0,j}^{(2,k)}$ enables to discriminate
multifractality from self-similarity.
Similar tests involving $\bar{\mu}_{i,F/A/U}^j$, $\bar{\mu}_{i,\nb_\ell}^j$ and
$\bar{\mu}_{i,\tb_r}^j$ in the definition of null hypotheses\reff{eq:StatsAssumptions} for $(i=1,2)$ were also performed.

\begin{figure}
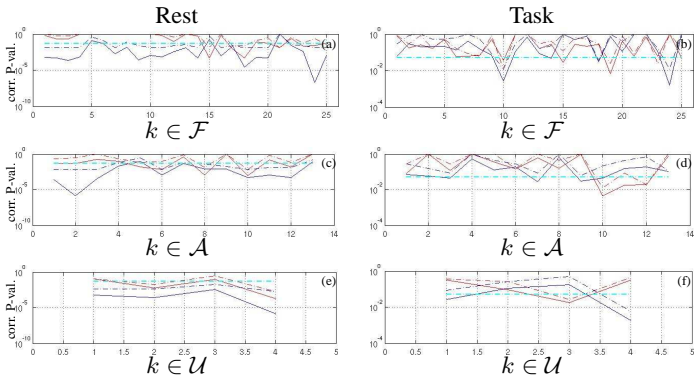

\bcc
\begin{tabular}{@{}c@{}cc}
\yaxis{{\tiny corr. P-val.}}& \figc[width=4.25cm]{figure5/spDict_corrPval_onesample_Rest_Fmaps2}&
\figc[width=4.25cm]{figure5/spDict_corrPval_onesample_TaskAD_Fmaps2}\\[-1.35cm]
&\hspace*{4cm}{\tiny (a)} & \hspace*{4cm}{\tiny (b)} \\[-.75cm]
 &{\small Rest} & {\small Task}\\[1.1cm]
&{\small $k\in\Fc$} & {\small $k\in\Fc$}\\
\yaxis{{\tiny corr. P-val.}}&\figc[width=4.25cm]{figure5/spDict_corrPval_onesample_Rest_Amaps2}&
\figc[width=4.25cm]{figure5/spDict_corrPval_onesample_TaskAD_Amaps2}\\[-1.35cm]
&\hspace*{4cm}{\tiny (c)} & \hspace*{4cm}{\tiny (d)} \\[.75cm]
&{\small $k\in\Ac$} & {\small $k\in\Ac$}\\
\yaxis{{\tiny corr. P-val.}}&\figc[width=4.25cm]{figure5/spDict_corrPval_onesample_Rest_Umaps2}&
\figc[width=4.25cm]{figure5/spDict_corrPval_onesample_TaskAD_Umaps2}\\[-1.35cm]
&\hspace*{4cm}{\tiny (e)} & \hspace*{4cm}{\tiny (f)} \\[.75cm]
&{\small $k\in\Uc$} & {\small $k\in\Uc$}
\end{tabular}\vspace*{-.75cm}
\ecc
\caption{Corrected p-values associated with one-sample Student-t~({\color{blue}{--}},{\color{red}{--}})
and WSR~({\color{blue}{-.}},{\color{red}{-.}}) tests performed for testing
\tiny$H_{0,j}^{(1,\cdot)}$\footnotesize~(\textcolor{blue}{blue} curves) and
\tiny$H_{0,j}^{(2,\cdot)}$\footnotesize~(\textcolor{red}{red} curves) on
the F-~(top), A-~(center) and U-maps~(bottom), respectively, for $j={\rm R}$~(left) and
$j={\rm T}$~(right). Significance level~($\alpha=.05$) is shown in \textcolor{cyan}{-~-}.
\label{fig:Pval_onesample_maps}}
\end{figure}

Analysis of statistical significance of F-maps regarding $H_{0,R}^{(1),k}$ showed that most components~(22/25) rejected this
null hypothesis at rest using T-test and thus were significantly self-similar~(see blue
curves in Fig.~\ref{fig:Pval_onesample_maps}(a)).
The task effect then induced a loss of significance in the vast majority of components as shown in
Fig.~\ref{fig:Pval_onesample_maps}(b): only four maps~($\fb_{10}$, $\fb_{14}$, $\fb_{18}$ and $\fb_{24}$) demonstrated
a significant level of self-similarity using T-test in the task-related dataset.
These maps are related to the motor, fronto-parietal and attentional~(parieto-temporal junction and IPS/FEF) networks.
Two out of them are lateralized in the left hemisphere. Statistical analysis of F-maps regarding $H_{0,R}^{(2),k}$
demonstrated that only six components~($\fb_{15}$, $\fb_{17}$, $\fb_{18}$, $\fb_{21}$, $\fb_{23}$, $\fb_{24}$)
rejected this null hypothesis at rest: see red curves in
Fig.~\ref{fig:Pval_onesample_maps}(a).
The task-related modulation tends to reduce the number of significant F-maps: As depicted in
Fig.~\ref{fig:Pval_onesample_maps}(b), only 3 components survived the T-test~($\fb_{10}$, $\fb_{15}$ and $\fb_{19}$)
in the task-related dataset.
Interestingly, $\fb_{10}$ and $\fb_{19}$ are likely to be involved in the auditory detection task
and the motor response since they belong to the Motor and Attentional networks.
Hence, a significant level
of multifractality is observed during task in components that were monofractal at rest. 
Besides, 
the level of multifractality remains significant in the ventral occipital cortex~($\fb_{15}$) irrespective of the
brain state and that a few components in the visual~($\fb_{17}$), fronto-parietal~($\fb_{18}$), temporal~($\fb_{21}$),
prefrontal~($\fb_{23}$) and attentional~($\fb_{24}$) networks became monofractal under the task effect.

Statistical analysis of A and U-maps regarding $H_{0,j}^{(1),k}$ showed the same behavior when switching from rest to task,
namely a strong decrease of the number of significant self-similar components~(from 10 to 4 and 4 to 2 for A/U-maps,
respectively): see blue curves in
Fig.~\ref{fig:Pval_onesample_maps}(c)-(d) and Fig.~\ref{fig:Pval_onesample_maps}(e)-(f), respectively.
Statistical analysis of A and U-maps regarding $H_{0,j}^{(2),k}$ also demonstrated a reduction of the number of
multifractal components in A/U-maps. Two artifactual components~($\ab_{10}$ and $\ab_{12}$) located in the white matter
remained consistently multifractal in both datasets and one undefined component~($\ub_3$) became significantly multifractal 
when switching from rest to task.
In all cases, a loss of significance is observed using WSR tests~(dash dotted curves) instead of T-tests~(solid curves)
indicating that there is no outlier in this group and thus that the Gaussian distribution hypothesis is
tenable.
    
\begin{figure}
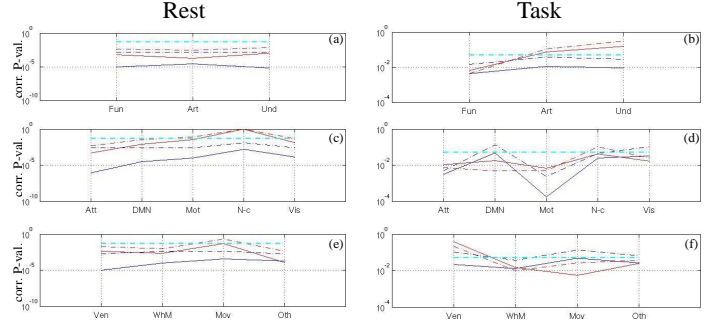

\bcc
\begin{tabular}{@{}c@{}cc}
\yaxis{{\tiny corr. P-val.}}& \figc[width=4.25cm]{figure6/spDict_corrPval_onesample_Rest_FAU2}&
\figc[width=4.25cm]{figure6/spDict_corrPval_onesample_TaskAD_FAU2}\\[-1.35cm]
&\hspace*{4cm}{\tiny (a)} & \hspace*{4cm}{\tiny (b)} \\[-.75cm]
&{\small Rest} & {\small Task}\\[1.25cm]
\yaxis{{\tiny corr. P-val.}}&\figc[width=4.25cm]{figure6/spDict_corrPval_onesample_Rest_Nets2}&
\figc[width=4.25cm]{figure6/spDict_corrPval_onesample_TaskAD_Nets2}\\[-1.35cm]
&\hspace*{4cm}{\tiny (c)} & \hspace*{4cm}{\tiny (d)} \\[1cm]
\yaxis{{\tiny corr. P-val.}}&\figc[width=4.25cm]{figure6/spDict_corrPval_onesample_Rest_Arts2}&
\figc[width=4.25cm]{figure6/spDict_corrPval_onesample_TaskAD_Arts2}\\[-1.35cm]
&\hspace*{4cm}{\tiny (e)} & \hspace*{4cm}{\tiny (f)} \\[1cm]
\end{tabular}\vspace*{-.75cm}
\ecc
\caption{Corrected p-values associated with one-sample Student-t~({\color{blue}{--}},{\color{red}{--}})
and WSR~({\color{blue}{-.}},{\color{red}{-.}}) tests performed for testing
\tiny$\bar{H}_{0,j}^{(1,\cdot)}$\footnotesize~(\textcolor{blue}{blue} curves) and
\tiny$\bar{H}_{0,j}^{(2,\cdot)}$\footnotesize~(\textcolor{red}{red} curves) on the
the averaged map types~(top), networks~(center) and artifact types~(bottom), respectively for $j={\rm R}$~(left) and
$j={\rm T}$~(right). Significance level~($\alpha=.05$) is shown in \textcolor{cyan}{-~-}.
\label{fig:Pval_onesample_Nets}}
\end{figure}

Then, we focused on the statistical analysis at different macroscopic scales, first by averaging all F/A and U-maps respectively
so as to derive a mean behavior for F/A/U-maps. Finally, we looked at functional networks
and artifact types in more details.
Blue curves in Fig.~\ref{fig:Pval_onesample_Nets}(a)-(b) report such results for the rest and task-related datasets,
respectively.
We still observed a significant level of self-similarity in all averaged groups~(blue curves) irrespective of
the brain state: $\bar{H}_{0,j}^{(1,F/A/U)}$ is systematically rejected for $j=({\rm R,T})$.
However, we still noticed a reduction of statistical significance induced by task irrespective of the map type.
More interestingly, we found at this macroscopic level that all averaged maps were multifractal at rest whereas
only the functional one remained multifractal during task: see red curves in
Fig.~\ref{fig:Pval_onesample_Nets}(a)-(b). Further, statistical analysis of functional networks defined in
Tab.~\ref{tab:table_MSDL_maps} was conducted to understand which network drives this effect.
When comparing p-values in Fig.~\ref{fig:Pval_onesample_Nets}(c)-(d) on functional networks, we observed that all
remained significantly self-similar in both states, while the DMN is close to the significance level
$\alpha\!=\!0.05$ during task~(blue curves). Regarding multifractality, only the non-cortical regions appeared monofractral at
rest and all networks kept a significant amount of multifractality during task. In contrast, this observation did
not hold for artifacts: when looking at Fig.~\ref{fig:Pval_onesample_Nets}(e)-(f) in detail, the signal related
to ventricles became monofractal during task.

\subsubsection{2-way repeated measures ANOVA}

In order to assess any significant change of self-similarity or multifractality between rest and task, 
we entered the subject-dependent parameter estimates $(\wh{c}_{i,k}^{j,s})$ in several 2-way repeated measures ANOVAs
involving two factors: brain state~(two values: $j={\rm R, T}$) and map type~(with varying number of values).
These ANOVAs were conducted separately for assessing self-similarity~($i=1$) and
multifractality~($i=2$) changes.
First six ANOVAs~(three for each parameter) were carried out by considering the F/A/U-maps as the second factor, respectively. This second factor
thus took a number of values that depends on the set under study: $\Fc$, $\Ac$ or $\Uc$. 
Results are summarized in Tab.~\ref{tab:ANOVA_local}. 
Regarding the analysis of self-similarity~($\wh{c}_{1,k}^{j,s}$ parameters),
a significant brain state effect appeared in all F/A/U-maps, and a significant map effect in the F and A-sets.
Significant interactions were found for the F and U-maps. 
This confirms that the level of self-similarity is not sufficient to disentangle functional networks from artifactual or
undefined maps.\\
As regards ANOVAs based on $\wh{c}_{2,k}^{j,s}$ parameters, a significant interaction for F-maps is found, thus indicating
that the averaged change in multifractality between rest and task is significant for functional maps only.
In summary, only F-maps exhibited significant interactions for both multifractal attributes.

\begin{table}
\bcc
\caption{2-way repeated measures ANOVA results based on the $\wh{c}_{i,k}^{j,s}$ parameters
for $i=\acc{1,2}$, $j=({\rm R, T})$, $s=1:S$ and $k\in\Fc$~(top), $k\in\Ac$~(middle), $k\in\Uc$~(bottom).
\label{tab:ANOVA_local}}
{\small 
\begin{tabular}{c|c|c|c|c}\hline\hline
{\small Level} & {\small Param.} &  {\small Source}  & {\small F score} &{\small p-val.}\\\hline
\multirow{3}{*}{\small F-maps} & \multirow{3}{*}{$\wh{c}_{1,k}^{j,s}$} & {\small State} &  {\small 9.54} &   {\small\bf 0.01}\\
&      &  {\small Map}   &  4.31   & {\bf 1e-09}\\
&  & {\small State} $\times$ {\small Map}  & 1.76 &   {\bf 0.02}\\\hline
\multirow{3}{*}{\small F-maps} & \multirow{3}{*}{$\wh{c}_{2,k}^{j,s}$} & {\small State}  & 0.13 &   0.73\\
&      &  {\small Map}    &   1.19  & 0.25\\ 
&  & {\small State} $\times$ {\small Map}  &   1.56  &  {\bf 0.04}\\\hline \hline
\multirow{3}{*}{\small A-maps} & \multirow{3}{*}{$\wh{c}_{1,k}^{j,s}$} & {\small State} &  5.73 &  {\bf 0.03}\\
&      &  {\small Map}   & 2.4  &  {\bf 0.008} \\
&  & {\small State} $\times$ {\small Map} &  1.32 &  0.21\\ \hline
\multirow{3}{*}{\small A-maps} & \multirow{3}{*}{$\wh{c}_{2,k}^{j,s}$} & {\small State} & 0.09 &  0.77\\
&      &  {\small Map}  &  2.4  &  {\bf 0.007}\\
&  & {\small State} $\times$ {\small Map} &  0.71 &  0.74 \\\hline\hline
\multirow{3}{*}{\small U-maps} & \multirow{3}{*}{$\wh{c}_{1,k}^{j,s}$} & {\small State}  & 5.39 &  {\bf 0.04}\\ 
&      &  {\small Map}  &    2.91  & 0.06\\
&  & {\small State} $\times$ {\small Map}   & 3.16 &  {\bf 0.04} \\\hline
\multirow{3}{*}{\small U-maps} & \multirow{3}{*}{$\wh{c}_{2,k}^{j,s}$} & {\small State} & 2.43e-05 &   0.99\\
&      &  {\small Map}  &  0.68 &  0.57\\
&  & {\small State} $\times$ {\small Map}   & 0.63 &  0.6\\\hline\hline
\end{tabular}
}
\ecc
\end{table}

Akin to the one-sample analyses above, we looked at a larger spatial scale, the functional network and artifact
type levels and performed similar ANOVAs, corresponding results are reported in Tab.~\ref{tab:ANOVA_global}.
While both functional networks and artifacts demonstrate a significant change in the
self-similarity parameter between rest and task, only functional networks made the map-type effect
significant. More importantly, the key feature for discriminating functional networks from artifacts relied on ANOVAs
based on $\wh{c}_{2,k}^{j,s}$ parameters. Indeed, a significant network effect and more importantly
a significant interaction between rest and task are observed in functional networks.

\begin{table}
\bcc
\caption{2-way repeated measures ANOVA results based on the $\wh{c}_{i,k}^{j,s}$ parameters
for $i=\acc{1,2}$, $j=({\rm R, T})$, $s=1:S$ and $k\in\Nc$~(top) and $k\in\Tc$~(bottom).
\label{tab:ANOVA_global}}
{\small 
\begin{tabular}{c|c|c|c|c}\hline\hline
{\small Level} & {\small Param.} &  {\small Source}  & {\small F score} &{\small p-val.}\\\hline
\multirow{3}{*}{\small Networks} & \multirow{3}{*}{$\wh{c}_{1,k}^{j,s}$} & {\small State}   & 9.78 &  {\bf 0.01}\\
     &       &  {\small Network}   & 4.18   & {\bf 0.006}\\
     &  & {\small State} $\times$ {\small Network}  & 1.09   & 0.37\\\hline
     \multirow{3}{*}{\small Networks} & \multirow{3}{*}{$\wh{c}_{2,k}^{j,s}$} & {\small State}   & 1.013  & 0.34 \\
     &       &  {\small Network}    & 3.18 &  {\bf 0.02}\\
     &       &  {\small State} $\times$ {\small Network} & 2.97 &  {\bf 0.03}\\\hline\hline
\multirow{3}{*}{\small artifacts} & \multirow{3}{*}{$\wh{c}_{1,k}^{j,s}$} & {\small State} & 4.85 &  {\bf 0.05}\\
&       & {\small artifact}  &   2.33 &   0.09\\
&  & {\small State} $\times$ {\small artifact}   &1.16   & 0.34\\\hline
\multirow{3}{*}{\small artifacts} & \multirow{3}{*}{$\wh{c}_{2,k}^{j,s}$} & {\small State} &     0.31  &  0.59\\
&       & {\small artifact}    &   1.03  &  0.39\\
&  & {\small State} $\times$ {\small artifact}  &  1.085&  0.37\\\hline\hline
\end{tabular}
}
\ecc
\end{table}

\subsubsection{Two-sample statistical tests}

\begin{figure}
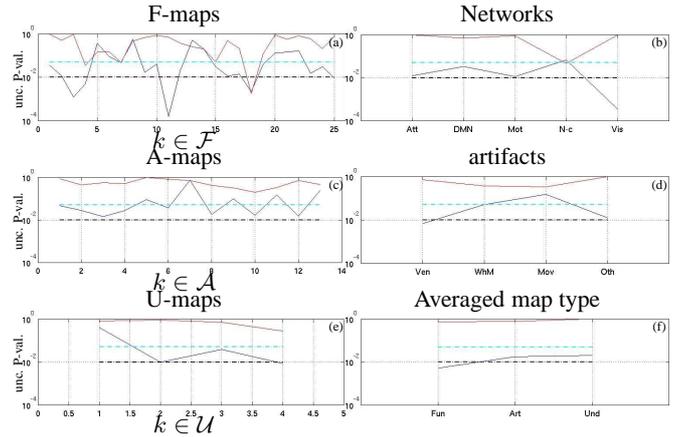

\bcc
\begin{tabular}{@{}c@{}c@{}c}
\yaxis{{\tiny unc. P-val.}}& 
\figc[width=4.25cm]{figure7/spDict_uncorrPval_twosample_Fmaps_2thresh}&
\figc[width=4.25cm]{figure7/spDict_uncorrPval_twosample_Nets_2thresh}\\[-1.55cm]
& \hspace*{4cm}{\tiny (a)} &\hspace*{4cm}{\tiny (b)}\\[-.75cm]
&{\small F-maps} & {\small Networks} \\[1.25cm]
&{\small $k\in\Fc$} & \\[.15cm]
\yaxis{{\tiny unc. P-val.}}& \figc[width=4.25cm]{figure7/spDict_uncorrPval_twosample_Amaps_2thresh}&
\figc[width=4.25cm]{figure7/spDict_uncorrPval_twosample_Arts_2thresh}\\[-1.55cm]
& \hspace*{4cm}{\tiny (c)} &\hspace*{4cm}{\tiny (d)}\\[-.75cm]
&{\small A-maps} & {\small artifacts} \\[1.25cm]
& {\small $k\in\Ac$} & \\[.15cm]
\yaxis{{\tiny unc. P-val.}}& 
\figc[width=4.25cm]{figure7/spDict_uncorrPval_twosample_Umaps_2thresh}&
\figc[width=4.25cm]{figure7/spDict_uncorrPval_twosample_FAU_2thresh}\\[-1.55cm]
& \hspace*{4cm}{\tiny (e)} &\hspace*{4cm}{\tiny (f)}\\[-.75cm]
&{\small U-maps} & {\small Averaged map type} \\[1.25cm]
& {\small $k\in\Uc$} & 
\end{tabular}\vspace*{-.75cm}
\ecc
\caption{Uncorrected p-values associated with two-samples Student-t test
performed for testing
\tiny$\wt{H}_{0}^{(1,\cdot)}$\footnotesize~(\textcolor{blue}{blue} curves) and
\tiny$\wt{H}_{0}^{(2,\cdot)}$\footnotesize~(\textcolor{red}{red} curves) on the
the F/A/U-maps~(left) and networks, artifacts and map types (right), respectively.
Significance levels~($\alpha_1=.01$ and $\alpha_2=.05$) are shown in {\bf \textcolor{black}{-~-}}
and \textcolor{cyan}{-~-}, respectively.
\label{fig:Pval_two_samples_maps}}
\end{figure}

To localize which maps are responsible for statistically significant ANOVA results, we finally performed
two-sample T-tests 
in which we tested the following null hypotheses:
\begin{equation}
\left\{
\begin{array}{ll}
\wt{H}_{0}^{(1,k)}: \mu_{1,k}^{\rm R}= \mu_{1,k}^{\rm T}, \; \forall k\in\Fc\cup\Ac\cup\Uc\\
\wt{H}_{0}^{(2,k)}:\mu_{2,k}^{\rm R}= \mu_{2,k}^{\rm T}, \; \forall k\in\Fc\cup\Ac\cup\Uc.
\end{array}
\right.
\label{eq:TwoSamples}
\end{equation}
We also conducted similar tests at the macroscopic level~($k\in\Nc\cup\Tc$) by replacing $\mu_{i,k}^j$ with $\bar{\mu}_{i,k}^j$ in the
null hypotheses~\eqref{eq:TwoSamples}.
The fluctuations in self-similarity being systematically in the same direction between rest
and task, we performed one-sided tests as regards the $\mu_{1,k}^j$'s while two-sided tests were
considered for the $\mu_{2,k}^j$'s: task-related positive and negative fluctuations of $\mu_{2,k}^j$ were actually observed
in Subsection~\ref{subsec:mean_effects}.
Fig.~\ref{fig:Pval_two_samples_maps}(a)-(b) shows the uncorrected p-values for the F-maps and networks,
respectively. 
We rejected $\wt{H}_{0}^{(1,k)}$ for $\pth{\fb_3,\fb_4,\fb_{11},\fb_{18},\fb_{25}}$ at a significance level set to
$\alpha_1\!=\!0.01$
and $\wt{H}_{0}^{(2,k)}$ for $\pth{\fb_4,\fb_7,\fb_{18}}$ at $\alpha_2\!=\!0.05$.
These components clearly explain significant results reported in Tab.~\ref{tab:ANOVA_local} about the changes in
self-similarity and multifractality that occurred in F-maps. Interestingly, among the latter, the
null hypothesis was rejected because of a large increase of multifractality in $(\fb_4,\fb_7)$.
In contrast, a decrease of multifractality was responsible for the rejection of $\wt{H}_{0}^{(2,k)}$ in $\fb_{18}$.
When setting $\alpha_2\!=\!\alpha_1\!=\!0.01$, only $\fb_{18}$ survived this threshold and thus remained the single functional
component for which a significant difference of self-similarity and multifractality was found between rest and task.
This component clearly drove the significant interaction reported in Tab.~\ref{tab:ANOVA_local}
for the change in multifractality in F-maps.
Fig.~\ref{fig:Pval_two_samples_maps}(b) also showed that the state effect reported in Tab.~\ref{tab:ANOVA_global}
on $(\wh{c}_{1,k}^{j,s})$ at the network level was driven by the attentional, motor and visual systems.
Last, the significant interaction reported in Tab.~\ref{tab:ANOVA_global}
on $(\wh{c}_{2,k}^{j,s})$ is explained by the non-cortical regions as shown in Fig.~\ref{fig:Pval_two_samples_maps}(b).

Fig.~\ref{fig:Pval_two_samples_maps}(c)-(d) shows the localization of the state effects reported in
Tabs.~\ref{tab:ANOVA_local}-\ref{tab:ANOVA_global} for the changes in self-similarity that occurred in artifacts at
the local and global levels. No A-map enabled to reject $\wt{H}_{0}^{(1,k)}$ at the $\alpha_1$ significance level but a
majority of A-maps~$\pth{\ab_{1:4},\ab_6,\ab_8,\ab_{10},\ab_{12}}$ contributed to the significant state effect observed
in Tab.~\ref{tab:ANOVA_local}. At the global artifact level, the ventricles appear as the main source of the significant
state effect reported in Tab.\ref{tab:ANOVA_global} for the change in self-similarity.
Also, no significant difference in multifractality was reported for artifacts whatever the observation
level~(A-maps or averaged artifacts).
Similarly, Fig.~\ref{fig:Pval_two_samples_maps}(e) enables us to show that $\ub_2$
and $\ub_4$ were the main sources of the significant state effect and interactions reported in Tab.~\ref{tab:ANOVA_local}
for the change in self-similarity. At the macroscopic level, we finally observed in Fig.~\ref{fig:Pval_two_samples_maps}(f)
that only the grand mean of functional maps leads to a significant modulation
of self-similarity between rest and task at level $\alpha_{1}$.

\section{Discussion}
\label{sec:discussion}

\subsection{Results interpretation}

This study analyzed in depth the scale-free properties of fMRI signals, using multifractal methodologies, and their modulations during rest and task both in functional networks and artifactual regions. The underlying goal was to finely characterize which properties are specific to functional
networks and which modulation can be expected for these networks from task-related activity.
Previous attempts in the literature~\cite{Cordes01,Leopold03,He11} focused on functional
networks without comparing results with the behavior of artifacts. The main reason comes from the fact that
seed region analyses were only conducted in such studies. Hence, no comparison with vascular or
ventricles-related signals was undertaken. 

Our results confirmed that fMRI signals are stationary and self-similar but not specifically in functional networks.
Also we showed that the amount of self-similarity significantly varies between rest and task not only in functional
networks involved in our auditory detection task with a motor response~(Attentional, Motor) but also quite
surprisingly in the visual
system and in some artifacts~(ventricles) and undefined maps. This observation led us to investigate the scale-free
structure of fMRI signals using richer models, namely multifractal processes, to which the WLMF toolbox is dedicated.
Our statistical results demonstrate first that fMRI signals are multifractal, second that interactions
between brain state and maps only occurred in F-maps and functional networks and third, that specific F-maps
such as in non-cortical regions demonstrated a statistically significant fluctuation between rest and task. 
This result shows that the concept of multifractality
permits to disentangle functional components from artifactual ones, in a robust and significant manner.

However, in contrast to self-similarity that systematically decreases with evoked activity,
multifractality decreases in cortical~($\fb_{18}$) but increases in non-cortical~($\fb_4,\fb_7$). 
Thus, task-related activity has no systematic impact with respect to increase/decrease of multifractality. 
Interestingly, we found a statistically non-significant trend towards a decrease of multifractality in
regions primarily involved in the task~($\fb_{12}$, $\fb_{24}$, $\fb_{25}$). However, the group size of this
study remains small~(12 subjects only) to achieve significant results, mainly because of the between-subject
variability and of the difficulty in estimating $\wh{c}_{2,k}^{j,s}$ parameters on short time series.

Further investigations beyond the scope of this paper are necessary to
find out any general trend on the direction change of multifractality with evoked-activity
by cross-correlating multifractal parameters with task-related activity~(e.g. group-level Z-scores)
and task performance. However, to derive reliable results for multifractality, 
a larger group of individuals will be considered and
a larger number of scans will be acquired while maintaining the same scanning time:
To this end, accelerated SENSE imaging
will be used together with recent reconstruction algorithms so as to improve temporal
resolution~\cite{Chaari10}.

\subsection{Monofractal scale-free EEG microstate sequences vs multifractal dynamics for RSN}

The results obtained in this contribution shows multifractal temporal dynamics in fMRI signals and thus
naturally lead to question the potential origins and generative mechanisms for this departure from the more
traditional longe range correlation modeling of scale invariance. 
A natural track to inspect consists of that of the relations between hemodynamic~(fMRI) and
electrical~(EEG) signatures for brain activity at rest.
This question has been intensively studied over the last
decade~\cite{Laufs03,Mantini07,Britz10,VandeVille10,Musso10,Yuan12},
first by measuring  cross correlations between fMRI data at rest
and EEG-informed regressors derived from the convolution of the EEG power signal in five well-identified frequency bands
($\delta \in (1,4)$~Hz, $\theta\in(4,7)$~Hz, $\alpha\in (8,12)$~Hz,
$\beta\in(13,30)$~Hz and $\gamma>30$~Hz) with the canonical HRF.
This approach revealed the negative correlation of $\alpha$-band activity with the attentional network
and the positive correlation with $\beta_2$-band with the default mode network~(precuneus and posterior
cingulate cortex)~\cite{Laufs03}. Also,~\cite{Mantini07}
showed that functional resting state networks have different EEG signatures which are not specific
to a given frequency band but are rather spread over several oscillations regimes~(e.g., correlation
between $\alpha$ and $\beta$ power in specific RSN), a consequence of the so-called oscillation
hierarchy~\cite{Buzsaki06} and the of phase-amplitude cross-frequency
coupling~\cite{He10}.
However, none of these works enable to explain the low frequency fluctuations~($<0.1$~Hz) or scale-free
dynamics of the fMRI signal at rest, because this phenomenon is much more widespread than oscillations.

Scale-free dynamics of brain electrical activity at rest has recently been
studied~\cite{VandeVille10} but not directly on raw data. Instead, EEG microscates that correspond to short periods~(100 ms)
during which the EEG scalp topography remains quasi-stable, have been first segmented. Remarkably, it
has been shown that only four different EEG microstate patterns are necessary to describe the ongoing
electrical brain activity at rest~\cite{Britz10} and that these four microscates correlate
with well-known RSNs, which were classically identified from fMRI dataset alone using group-level ICA.
This demonstrated that the EEG microstate with rapid fluctuations might be considered as the
electrophysiological signature of intrinsic functional connectivity patterns. The investigation of
scale-free dynamics was thus performed on the EEG microstate sequence to understand how fast the microscates
are changing and what kind of correlation structure~(short or long range) they
bring~\cite{VandeVille10}.

The recent finding that EEG microscate sequences reveal purely monofractal
dynamics~\cite{VandeVille10}, irrespective of the data filtering, may lead to conclude that the same monofractal 
behaviour in the fMRI signature of RSN (strongly correlate with these microstates) should be expected, if one
assumes a linear and time invariant HRF model for the neurovascular coupling. 
However, the results obtained in the present contribution can be considered not only as evidence in favor
of multifractality in fMRI data, but also as evidence that this multifractal effect is discriminant of cortical
versus non cortical regions and characteristic of functional network with respect to modulation under task.

Several factors may explain this apparent discrepancy.
First, an accurate comparison of both sets of result would require a precise match of the range of
scales~(or frequencies) within which scale invariance is analyzed and corresponding parameters measured. 
Here, the selected range of frequencies
corresponds to ([.008, .063]Hz), while the monofractal behavior of EEG microstate sequences was exhibited
on a distinct frequency range ie. ([.063, 3.9]Hz) in~\cite{VandeVille10}.
Comparison of scaling properties requires that the same frequency range is selected but
this constraint is clearly not tenable across modalities like EEG and fMRI given the fMRI sampling rate.

Second, it is indeed very unlikely that a linear and time invariant filtering may create multifractality in
fMRI starting from a monofractal electrophysiological signal in EEG. 
The general issue of the relations between~(linear and non linear)
filtering and multifractality were barely studied theoretically so far but 
interestingly,~\cite{ABRY:2012:A} has shown that simple nonlinear filter can turn
mono- into multifractality.
Hence, another putative origin for the apparent contradiction between our findings and
those in~\cite{VandeVille10} lies in refined descriptions of HRF model
by nonlinear dynamical systems~(e.g., Balloon model)~\cite{Buxton98,Buxton04}.
Of course, linear and stationary approximations like the canonical HRF model~\cite{Glover99} or nonparametric
alternatives~\cite{Vincent10,Chaari11}
have been validated but only on evoked activity and considering inter-stimulus intervals larger than 3~s.
For shorter ISIs, nonlinear hemodynamics has turned out to be a valid property~\cite{Liu00}.
In this context, habituation or repetition supression effects
may occur and induce a sublinear hemodynamic response, which would modify scaling
properties~\cite{Dehaene-Lambertz06,Ciuciu09}. Hence, by modelling the sequence of EEG
transient brain states as a series of short \emph{time epochs}, this could induce nonlinearities in the
hemodynamic system that could explain the switch from purely fractal EEG microstates to multifractal
signatures in the corresponding RSNs.

Third, instead of segregating EEG microstates in multiple groups
based upon the maximal spatial dissimilarity between groups~\cite{Britz10,Musso10},
a more recent analysis of joint EEG/fMRI resting state data has revealed
a larger number~(thirteen) of EEG microstates that show temporal independence from each
other~\cite{Yuan12}. In this latter work, all resting state networks including visual,
motor, auditory, attention, saliency and default mode networks were characterized by a specific
electrophysiological signature involving several EEG microstates. This clearly indicates that the original
analysis of scale-free dynamics for EEG microscates done in~\cite{VandeVille10} should be revisited on
this larger number of metastable states to disentangle whether multifractality in this larger set of
microstates has been discarded due to averaging effects. 
It is actually clear that the sequence mixing thirteen different microstates may generate richer
\emph{singularities}~(abrupt changes between microstates) than the ones relying on four microstates only.
Fourth, the temporal signatures of EEG microstates found in~\cite{VandeVille10,Musso10}
are correlated in time since the spatial similarity was the key factor to identify them. As a consequence,
the microstate sequences is correlated too and might loose some singularities
that could be found out in the microstate sequences generated by~\cite{Yuan12}.
Finally, the presence of multifractality in resting state~(and task-related) MEG
data has been evidenced in the sensor space in~\cite{Zilber12}.
These findings open new research avenues: For instance, it is natural to explore whether the observed multifractal properties can be related multiplicative cascade processes, that is to one one of the only practical mechanism known to generate multifractal dynamics, or to investigate whether this cascade takes place at
meso or macroscopic scales, as well as to figure out how brain networks could implement such cascade mechanisms. 
This topic is beyond the scope of the present contribution,  however the log-normal statistics of neuronal firing
rate could provide us with a first clue to uncover any generative process underlying multifractal dynamics.

\subsection{Stationarity vs non-stationarity of the RSN dynamics}

Recent results in resting state fMRI reveal temporally independent functional modes of spontaneous brain
activity~\cite{Smith12}
and postulate the presence of temporally non-stationary modes in part of the default mode network
by resorting to high temporal resolution fMRI. 
While stationarity receives a unique and clear definition, non stationarity can correspond to a bunch of
different situations; for example, non-stationarity might \emph{(i)} refer to an apparent change over time
in the correlation between two regions or \emph{(ii)} refer to changes in the mean and/or variance in the time
course of a functional network.

The wavelet based analysis of scaling proposed here already addresses a number of
such situations. The fact that the estimated Hurst coefficient of fMRI time series remains consistently below 1
indicates that fMRI signals at hand here are better modeled as a stationary step process $Y$
rather than as a non stationary random walk $X$.
Further, wavelet analysis are known to bring robustness against various forms of non stationarities,
such as smooth trends superimposed to data, to mean or variance modulation~(cf {\em (ii)}).
The multifractal analysis performed here is thus not impaired by such form of non stationarities.
This leaves open issues such as the presence of oscillations superimposed to scaling.
Given that time series are very short, the use of formal stationarity test will lack power and are not
likely to reject stationarity.
Further, in all the analysis conducted in the present work, no evidence of non-stationarity in the
fMRI time series at hand were evidenced. This is in agreement with what has been reported in~\cite{He11}
in an fMRI ROI-based analysis. Finally, previous attempts to scale-free analysis of densely sampled fMRI datasets in
time~(using the EVI sequence~\cite{Rabrait08} already confirmed the validity of a the stationarity assumption;
see~\cite{Ciuciu08}.

\section{Conclusion}\label{sec:conclusion}

We uncovered multifractal scale-free dynamics of fMRI time series over four octaves~(15s.-125s.)
both in functional networks and in artifacts. We then disentangled
functional components from artifactual ones in a robust and significant manner by demonstrating
that only the former gave rise to significant modulations of the multifractal attributes between rest and
task-related activity. Variability in human performance scores also generally exhibits power law distributions,
whose strength (or exponent) is often modulated across conditions and tasks~\cite{Holden11}.
This paves the way towards future works devoted to investigating the extent to which behavioral properties
are correlated with the change of scale-free dynamics in neuroimaging time series~(MEG, fMRI) acquired during
multisensory learning~\cite{Seitz07}.

\section{Acknowledgements}\label{sec:acknowledge}

The authors thank the French National Research Agency~(ANR) for its financial support to the SCHUBERT~(ANR-09-JCJC-0071)
young researcher project~(2009-13). The authors are also grateful to the anonymous reviewers and the associate editor, Dr.
BJ He, for their remarks and criticisms that helped us to improve the manuscript and enlarge the readership.
Finally, we would like to warmly thank Dr. V. van Wassenhove for her careful rereading and constructive comments.

{\small
\bibliographystyle{model1b-num-names}
\bibliography{finalVersion_Ciuciu_FrontiersFractalPhysioRT_with_macros}
}

\end{document}